\documentclass[11pt,leqno,twoside]{article}
\def\cvd{\hfill$\Box$}

\def\corr{\longleftrightarrow}
\def\w{\dot{w}}
\def\v{{\rm v}}

\def\int{\mathbb{Z}}
\def\C{\mathbb{C}}
\def\R{\mathbb{R}}

\def\Z{\mathbb{Z}}
\def\N{\mathbb{N}}
\def\e{\varepsilon}
\def\Ue{{\cal U}_{\varepsilon}({\mathfrak g})}

\def\im{\text{Im}}

\def\Frac{\text{Frac}~}
\def\O{{\cal O}}

\def\h{{\mathfrak h}}

\def\Class{{\noindent\bf Class }}

\def\Classes{{\noindent\bf Classes }}

\def\proof{\noindent{\bf Proof. }}
\def\Pf{\proof}

\def\pf{\proof}

\def\rk{{\rm rk}}
\def\ov{\overline}
\def\bG{\overline G}

\def\a{\alpha}

\def\b{\beta}
\def\d{\delta}
\def\l{\lambda}
\def\g{\gamma}
\def\s{\sigma}

\def\<#1{\langle #1\rangle}
\def\o#1{|\ \! #1\!\ |}
\def\Ra{\Rightarrow}
\def\wJ{w_{\!_J}}

\def\vuoto{\varnothing}

\usepackage{times}
\usepackage{graphics}
\usepackage{amssymb}
\usepackage{amscd}
\usepackage{amsmath}
\usepackage{fancyhdr}

\input xy
\xyoption{all}

\title{On the coordinate ring of spherical conjugacy classes}

\newtheorem{theorem}{Theorem}[section]
\newtheorem{lemma}[theorem]{Lemma}
\newtheorem{corollary}[theorem]{Corollary}
\newtheorem{proposition}[theorem]{Proposition}
\newtheorem{definition}[theorem]{Definition}
\newtheorem{remark}[theorem]{Remark}

\newcounter{tigre}
\def\totable{\refstepcounter{tigre}\arabic{tigre}}

\fancypagestyle{plain}{%
\fancyhf{}

}

\author{Mauro Costantini\\
Dipartimento di Matematica Pura ed Applicata\\
Torre Archimede - via Trieste 63 - 35121 Padova - Italy\\
email: costantini@math.unipd.it }
\date{}
\begin{document}
\baselineskip=16pt
\maketitle
\begin{abstract}
  Let $G$ be a simple algebraic group over an
  algebraically closed field $k$ of characteristic zero and $\O$ be a spherical conjugacy class of $G$.  We determine the decomposition of the coordinate ring $k[\O]$ of $\O$ into simple $G$-modules.

  \end{abstract}
  %-----------------------------------------------------------------------------------
%\section*{Introduction}
\section{Introduction}

\newcounter{equat}
\def\theequat{(\arabic{equat})}
\def\equat{\refstepcounter{equat}$$~}
\def\endequat{\leqno{\boldsymbol{(\arabic{equat})}}~$$}

\newcommand{\elem}[1]{\stackrel{#1}{\longto}}
\newcommand{\map}[1]{\stackrel{#1}{\to}}
\def\imp{\Rightarrow}
\def\Imp{\Longrightarrow}
\def\iff{\Leftrightarrow}
\def\Iff{\Longleftrightarrow}
\def\to{\rightarrow}
\def\longto{\longrightarrow}
\def\injto{\hookrightarrow}
\def\rtordu{\rightsquigarrow}

In   \cite{CCC} we proved the De Concini-Kac-Procesi
conjecture on the quantized enveloping algebra ${\cal U}_\varepsilon(\frak g)$ (introduced
in \cite{DC-K}) for simple 
${\cal U}_\varepsilon(\frak g)$-modules over spherical conjugacy classes of $G$
(we recall that a conjugacy class $\O$ in $G$ is
called {\it spherical} if a Borel subgroup of $G$ has a
dense orbit in $\O$): our main tool was the representation theory of the quantized Borel
subalgebra $B_{\e}$ introduced in \cite{3}. 

To fix the notation, $G$ is a complex
simple simply-connected algebraic group, $\mathfrak g$ its Lie algebra, $B$ a Borel subgroup of $G$, $T$ a maximal torus of $B$, $B^-$ the Borel subgroup opposite to $B$, $\{\alpha_1,\ldots,\alpha_n\}$ the set of simple roots with respect to the choice of $(T,B)$.
Let $W$ be the Weyl group of $G$ and let us denote by
$s_i$ the reflection corresponding to the simple root $\alpha_i$: 
$\ell(w)$ is the length of the element $w\in W$ and $\rk(1-w)$ is the rank of $1-w$ in the geometric representation of $W$. 

The representation theory of ${\cal U}_\varepsilon(\frak g)$ is related to the stratification of $G$ given by conjugacy classes, while
the representation theory of $B_{\e}$ is related to the stratification $\{X_w\mid w\in W\}$ of $B^-$, where 
 $X_w=B^-\cap B{w}B$ for every $w\in W$ (each $X_w$ is an affine variety
of dimension $n+\ell(w)$). 
We proved that for every spherical conjugacy class $\O$ in $G$, there exists $w\in W$ such that $\O\cap X_w\not=\vuoto$ and $\ell(w)+rk (1-w)=\dim{\cal O}$: this then allows to prove the De Concini-Kac-Procesi
conjecture for simple ${\cal U}_\varepsilon(\frak g)$-modules over elements in $\O$. In fact we proved also a result in the opposite direction, giving therefore a characterization of spherical conjugacy classes in terms of the Weyl group (\cite{CCC}, Theorem 25): 

let $\O$ be a conjugacy class of $G$ and $w=w(\O)$ be the unique element in $W$ such that $\O\cap BwB$ is dense in $\O$. Then $\O$ is spherical if and only if $\dim{\cal O}=\ell(w)+rk (1-w)$. 

Moreover $w$ is always an involution (see \cite{CCC}, Remark 4, \cite{Car}, Theorem 2.7). From this result we conjectured that, for a spherical $\O$,  the decomposition of the ring $\C[\O]$ of regular functions on $\O$ (to which we refer as to {\it the coordinate ring} of $\O$) as a $G$-module should be strictly related to $w(\O)$. This is the motivation for the present paper.

We recall that $\C[\O]$ is multiplicity-free, so that in order to obtain the decomposition of $\C[\O]$ into simple components one has just to determine which simple modules occur in $\C[\O]$:
$$
\C[\O] \underset G\cong \bigoplus_{\l\in\l(\O)} V(\l)
$$
where for each dominant weight $\l$, $V(\l)$ is the simple $G$-module of highest weight $\l$ (if $\l\in\l(\O)$ we say that $\l$ {\em occurs} in $\C[\O]$).

The decomposition of the coordinate ring $\C[X]$ for $G$-varieties $X$ has been investigated by various authors. 
If $\l$ is a non-zero highest weight, and $v\in V(\l)$ is a non-zero highest weight vector, then $\C[G.v]$ is isomorphic to 
$\underset{n\geq 0}\oplus  V(n\l)^\ast$ (\cite{VP}, Theorem 2). In particular this determines $\C[\O]$ for the minimal unipotent orbit of $G$. For a unipotent class in $G$ (equivalently nilpotent orbit in $\frak g$) McGovern (\cite{MG1}, Theorem 3.1) decribes $\C[\O]$ in terms of induced building blocks from a certain Levi subgroup of $G$ (via sheaf cohomology on $G/Q$, $Q$ a parabolic subgroup of $G$ associated to $\O$): it is then possible to obtain multiplicities of simple $G$-modules in $\C[\O]$ as an alternating sum of certain partition functions.
In the same paper the author gives a formula for $\C[\hat \O]$, where $\hat \O$ is the simply-connected cover of $\O$ (\cite{MG1}, Theorem 4.1). Then in \cite{MG2} there are tables for the sets of simple modules in $\C[\hat \O]$ for spherical unipotent classes in the classical groups (and conjecturally in the exceptional groups). For type $F_4$ the monoid $\l(\O)$ has been described in \cite{Broer1} for all spherical unipotent classes. For the maximal spherical unipotent class $\O$ in $E_8$, it has been shown in \cite{vogan}, Theorem 1.1, that every simple $G$-module occurs in $\C[\O]$ (so that $\O$ is a model orbit). 
In \cite{pany5}, Panyushev gives tables for the sets of simple modules  for (spherical) nilpotent orbits of height 2 (and conjecturally for height 3).
In \cite{Luna} the author describes explicitly the structure of principal model homogeneous spaces. For semisimple spherical classes, the description of $\l(\O)$ may be deduced from the tables in \cite{Kramer}. See also \cite{Vust}, Th\'eor\`eme 3, where symmetric varieties are considered.

The
main result of this paper is the following:

\medskip
\noindent
{\bf Theorem.}\ {\it ~Assume $\O$ is a spherical conjugacy class in $G$, and let $w=w(\O)$. Then a dominant weight $\lambda$ occurs in $\C[\O]$ if and only if 
$w(\lambda)=-\lambda$ and $\lambda(S_\O)=1$}. 

\medskip
Here $S_\O$ is a certain (finite) elementary abelian 2-subgroup of $T$ which we determine for every spherical conjugacy class, describing therefore explicitly $\l(\O)$: see tables $1,\ldots,26$. In particular we completely solve the problem of determining the simple modules occurring in $\C[\O]$ for unipotent classes (\cite{Jan}, 8.13, Remark 2),  and obtain the decomposition of $\C[\O]$ for conjugacy classes of mixed elements.

\medskip

Our proof is based on the deformation result obtained by Brion in  \cite{deformation}. We have $\C[\O]=\C[G/H]=\C[G]^H$, where $H$ is the centralizer of an element of $\O$ in $G$. There exists a flat deformation of $G/H$ to a quotient $G/H_0$, where $H_0$ contains the unipotent radical $U^-$ of $B^-$. We determine the decomposition of $\C[G/H_0]$ into simple components (i.e. we determine $\l(G/H_0)$), relating the group $H_0$ with $H$ via the theory of elementary embeddings (\cite{LV}, \cite{BLV}). We then prove the crucial fact that $\l(\O)$ is saturated (\cite{pany4}, \S 1.3), so that $\C[G/H]=\C[G/H_0]$ as $G$-modules. We also determine the decomposition of the coordinate ring $\C[\hat\O]$ for the simply-connected cover $\hat\O$ of $\O$, and of $\C[\ov\O]$.

The paper is structured as follows. In Section \ref{bepi} we introduce the notation. In Section \ref{mth} we recall some basic
facts about spherical varieties and we prove the main theorem.
In Section \ref{vari} we determine the group $S_\O$ for the  spherical conjugacy classes in the various groups, determining therefore the monoid $\lambda(\O)$, and also $\l(\hat O)$. In Section \ref{non-normal} we consider the coordinate ring $\C[\ov\O]$ of the closure of $\O$. It is well known that $\C[\ov\O]=\C[\O]$ if and only if $\ov\O$ is normal: we list all cases in which the spherical conjugacy class $\O$ has normal closure and we determine $\l(\ov \O)$ for the classes with non-normal closure.
In section \ref{generale} we consider the case when $G$ in not necessarily simply-connected.

All the results and proofs of this article remain valid for
$G$ a simple simply-connected algebraic group over an
  algebraically closed field $k$ of characteristic zero. 

{\center\bf Acknowledgements.}
It is a pleasure to thank P.\ Bravi and M.\ Brion for helpful
discussions and suggestions. 

\section{Preliminaries}\label{bepi}

We denote by $\C$ the complex numbers, by $\R$ the reals, by
$\Z$ the integers and by $\N$ the natural numbers. 

Let $A = (a_{ij})$ be a finite indecomposable
Cartan matrix of rank $n$. 
To $A$ there is associated a root system $\Phi$, a simple Lie 
algebra
${\mathfrak g}$ and a simple simply-connected algebraic group $G$ over 
$\C$. We fix a maximal torus $T$ of $G$, and a Borel subgroup $B$
 containing $T$:  $B^-$ is the Borel subgroup opposite to $B$, $U$ (respectively $U^-$) is the unipotent radical of $B$ (respectively of $B^-$). If $\chi$ is a character of $T$, we still denote by $\chi$ the character of $B$ which extends $\chi$.
 We denote by $\h$ the Lie
algebra of $T$. Then  $\Phi$ is the set of roots relative
to $T$, and $B$ determines the set of positive roots  
$\Phi^+$, and the simple roots $\Delta=\{\alpha_1,\ldots,\alpha_n\}$.  We fix a total ordering on $\Phi^+$ compatible with the height function.
We shall use the numbering and the description of the simple roots in terms of the canonical basis $(e_1,\ldots, e_k)$ of an appropriate $\R^k$ as in \cite{bourbaki}, Planches I-IX. For the exceptional groups, we shall write $\b=(m_1,\ldots,m_n)$ for $\b=m_1\a_1+\ldots+m_n\a_n$.

If $\gamma$ is a character of $T$, we shall also denote by $\gamma$ the corresponding linear form $(d\gamma)_1$ on $\frak h$.
The real subspace of ${\h}^\ast$ spanned by the
roots  is a Euclidean space $E$, endowed with the scalar product
$(\a_i,\alpha_j) = d_ia_{ij}$. Here
$\{d_1,\ldots,d_n\}$  are
relatively prime positive integers such that if $D$ is the
diagonal matrix with entries $d_1,\ldots,d_n$, then $DA$ is
symmetric. 
$P$ is the weight lattice, $P^+$ the monoid of dominant weights and $W$ the Weyl
group; $s_i$ is the simple reflection associated to $\alpha_i$,
$\{\omega_1,\ldots,\omega_n\}$ are the fundamental weights, $w_0$ is the longest element of $W$. In the expression $\l=\sum_i k_in_i\omega_i$ we always assume $k_i$'s and $n_i$'s in $\N$.
 If $V$ is a $G$-module, $v\in V$, $f\in V^\ast$, then the matrix coefficient $c_{f,v}:G\to \C$ is defined by $c_{f,v}(g)=f(g.v)$ for $g\in G$.
We consider the action of $G\times G$ on $\C[G]$  
$$
((g,g_1).f)(c) = f(g^{-1}cg_1)
$$
for $c$, $g$, $g_1\in G$, $f\in \C[G]$. The algebraic version of the Peter-Weyl theorem gives the decomposition
\begin{equation}
 \C[G]= \bigoplus_{\lambda \in P^+} V(-w_0\lambda)^\ast \otimes
V(-w_0\lambda)
\label{decomposizione}
\end{equation}
We put
$\Pi=\{1,\ldots,n\}$ and we fix a Chevalley basis $\{h_{i},  i\in \Pi; e_\a,\a\in\Phi\}$
 of $\frak g$.
We shall denote by
$\check{\omega}_i$, for $i=1,\dots, n$, the elements in ${\mathfrak h}$
defined by 
$
\alpha_j(\check{\omega}_i)=\d_{ij}
$ (recall that 
$
\omega_j(h_i)=\d_{ij}
$) for $j=1,\ldots,n$.
As usual we put 
$
\<{x,y}=\frac{2(x,y)}{(y,y)}$. 

We use the notation $x_\a(k)$, $h_\a(z)$, for $\a\in \Phi$, $k\in \C$, $z\in\C^\ast$ as in \cite{yale}, \cite{Carter1}. For $\a\in \Phi$ we put $X_\a=\{x_\a(k)\mid k\in \C\}$, the root-subgroup corresponding to $\a$, and $H_\a=\{h_\a(z)\mid z\in \C^\ast\}$. For $h\in \frak h$ we put $H_h=\exp \C h$.
  We identify $W$ with $N/T$, where $N$ is the normalizer of $T$:  given an element $w\in W$ we shall denote
a representative of $w$ in $N$ by $\dot{w}$. We choose the $x_\a$'s  so that, for all $\a\in \Phi$, $n_\a=x_\a(1)x_{-\a}(-1)x_\a(1)$
lies in $N$ and has image the reflection $s_\a$ in $W$. Then 
\begin{equation}
x_\a(\xi)x_{-\a}(-\xi^{-1})x_\a(\xi)=h_{\a}(\xi)n_\a\quad,\quad
n_\a^2=h_{\a}(-1)
\label{relazioni}
\end{equation}
for every $\xi\in\C^\ast$, $\a\in\Phi$ (\cite{springer}, Proposition 11.2.1).

 We put $T^w=\{t\in T\mid wtw^{-1}=t\}$, $T_2=\{t\in T\mid t^2=1\}$. In particular $T^w=T_2$ if $w=w_0=-1$.

For algebraic groups we use the notation in \cite{Hu2}, \cite{Carter2}. In particular, 
for $J\subseteq \Pi$, $\Delta_J=\{\a_j\mid j\in J\}$, $\Phi_J$ is the corresponding root system, $W_J$ the Weyl group, $P_J$ the standard parabolic subgroup of $G$, $L_J=T\<{X_\a\mid \a\in \Phi_J}$ the standard Levi subgroup of $P_J$.
For $z\in W$ we put 
$
U_z=U\cap z^{-1}U^-z$.
Then the unipotent radical $R_uP_J$ of $P_J$ is $U_{w_0\wJ}$,
where $\wJ$ is the longest element of $W_J$.
Moreover
$
U\cap L_J=U_{\wJ}$
is a maximal unipotent subgroup of $L_J$.

If $\Psi$ is a subsystem of type $X_r$ of $\Phi$ and $H$ is the subgroup generated by $X_\a$, $\a\in \Psi$, we say that $H$ is a $X_r$-subgroup of $G$.

 If $X$ is an algebraic variety, we denote by $\C[X]$ the ring of regular functions on $X$. If $X$ is a multiplicity-free $G$-variety, then we denote by $\l(X)$ the set of dominant weights  occurring in $\C[X]$, i.e.  $\l\in P^+$ such that $\C[X]$ contains (a copy of) $V(\l)$. If $x\in X$ we denote by $G.x$ the $G$-orbit of $x$ and by  $G_x$ the isotropy subgroup of $x$ in $G$. If the homogeneous space $G/H$ is spherical, we say that $H$ is a spherical subgroup of $G$.

If $x$ is an element of a group $K$ and $H\leq K$, we shall also denote by $C(x)$ the centralizer of $x$ in $K$, and by $C_H(x)$ the centralizer of $x$ in $H$. If $x$, $y\in K$, then $x\sim y$ means that $x$, $y$ are conjugate in $K$.
For unipotent classes in exceptional groups we use the notation in \cite{Carter2}. We use the description of centralizers of involutions as in \cite{Iwa}.

\section{The main theorem}\label{mth}

Let $\O$ be a spherical conjugacy class. Our aim is to determine $\l(\O)$. For this purpose if $H$ is the centralizer of an element in $\O$, we have $\C[\O]=\C[G/H]=\C[G]^H$ and, from  $(\ref{decomposizione})$, 
$$
\C[G]^H=\bigoplus_{\lambda \in \l(\O)} V(-w_0\l)^\ast\otimes u_\l
$$
 where $0\not=u_\l\in V(-w_0\l)^H$ (\cite{PV}, Theorem 3.12). We start by considering in general a spherical homogeneous space $G/H$. Without loss of generality we may assume $BH$ dense in $G$. By \cite{deformation}, Theorem 1, there exists a 
(flat) deformation of $G/H$ to a homogeneuos (spherical) space $G/H_0$, where $H_0$ contains a maximal unipotent subgroup of $G$ (such an homogeneous space is called {\em horospherical}, and $H_0$ a horospherical contraction of $H$).  An {\em elementary embedding} of $G/H$ is a pair $(X,x)$ where $X$ is a normal algebraic $G$-variety, $x\in X$ is such that $G.x$ is dense in $X$, $G_x=H$ and $X\setminus G.x$ is a $G$-orbit of codimension 1 (\cite{BP}, 2.2). In \cite{deformation} Brion constructs a $G\times \C^\ast$-variety and a flat $G\times \C^\ast$-morphism $p:Z\to\C$ (where $G$ acts trivially on $\C$ and $\C^\ast$ acts via  homotheties) such that $p^{-1}(\C^\ast)\cong G/H\times \C^\ast$ and $p^{-1}(0)\cong G/H_0$ (\cite{deformation}, Theoreme 1, \cite{BP} \S 3.11). One may consider $Z$ as an elementary embedding $(Z,z)$ of $(G\times \C^\ast)/(H\times 1)$, with closed orbit $(G\times \C^\ast)/(H_0\times \C^\ast)$;   $H\times 1$ is the isotropy subgroup of $z$, $H_0\times \C^\ast$ is the isotropy subgroup of an element in the closed orbit (\cite{BP}, proof of Corollaire 3.7). Let $P=P_J$ be the parabolic subgroup {\em associated} to $H$,  $P=\{g\in G\mid gBH=BH\}$, and let $L$ be a Levi subgroup (which we may assume equal to $L_J$, by taking an appropriate  conjugate of $H$ instead of $H$) of $P$ {\em adapted} to $H$ (\cite{BP}, 2.9): in particular
\begin{equation}
P\cap H=L\cap H\quad,\quad
L'\leq H\label{adatto}
\end{equation}
Then $P\times \C^\ast$ is the parabolic subgroup of $G\times \C^\ast$ associated to $H\times 1$ and $L\times \C^\ast$ is a Levi subgroup adapted to $H\times 1$ (\cite{BP}, Corollaire 3.7 and its proof).

By \cite{BP}, Proposition 3.10, i), we have
$
H_0\times \C^\ast=(R_uQ\times 1)(L\times \C^\ast \cap H_0\times \C^\ast)
$ 
where $Q$ is the opposite parabolic subgroup of $P$ with respect to $L$, so that
\begin{equation}
H_0=(R_uQ)(L\cap H_0)
\label{prima}
\end{equation}

We show that $L\cap H=L\cap H_0$. Let $L=CL'$, where $C$ is the connected component of the centre of $L$. Then $L'$ is contained also in $H_0$, by \cite{BP}, Th\'eor\`eme 3.6.

By \cite{BP}, Proposition 3.4, $Z$ contains an open
$P \times {\C}^*$-stable subset isomorphic to $R_uP \times W$
where $W$ is $L \times {\C}^*$-stable and meets the closed orbit, and $(W,z)$ is an elementary embedding of the torus $(C\times \C^\ast)/(C\cap H\times 1)$ (\cite{BLV}, proof of Lemme 4.2). Then $f=p_{|W}:W\to \C$ is a
$(C\times \C^\ast)$-equivariant flat morphism such that $f^{-1}(\C^\ast)\cong C/C\cap H\times \C^\ast$ and $f^{-1}(0)\cong C/H_0\cap C$.
So the coordinate rings of these orbits are isomorphic $C$-modules
and it follows that 
the isotropy  
groups of all points of $W$ are the same.
In particular 
\begin{equation}
C\cap H=C\cap H_0\label{rigido}
\end{equation}

With the above notation we prove

\begin{theorem}\label{nane2} Let $H$ be a spherical subgroup of $G$ such that $BH$ is dense in $G$ and $L=L_J$ is a Levi subgroup adapted to $H$. Then 
$
H_0=R_uQ\,(L\cap H)=\<{U^-,U_{\wJ}, C\cap H}
$. 
\end{theorem}
\pf
By (\ref{rigido}) we have 
$$
L\cap H_0=L'C\cap H_0=L'(C\cap H_0)=L'(C\cap H)=L'C\cap H=L\cap H
$$
so that by (\ref{prima}) we conclude.
\cvd

\begin{definition}\label{sferica} We put $\tilde\l(G/H)=\l(G/H_0)$.
\end{definition}
Note that $\l(G/H)\leq \tilde\l(G/H)$ since $BH$ is dense in $G$, and more generally $\Z\,\l(G/H)\cap P^+\leq \tilde\l(G/H)$ (\cite{pany4}, part 2 of the proof of Proposition 1.5). 
Moreover
\begin{equation}
\l(G/H_0)=\{\l\in P^+\mid \l(T\cap H)=1\}
\label{questo}
\end{equation}
since $\prod_{j\in J} H_{\a_j}\leq H$ and $X_{\a_j}.v_{-\l}=v_{-\l}$ if $(\l,\a_j)=0$ (here $v_{-\l}$ is a lowest weight vector of weight $-\l$ in $V(-w_0\l)$).  Also
$B\cap H\leq P\cap H=L\cap H$, so that $B\cap H=U_{\wJ}(T\cap H)$. If $\l\in \tilde\l(G/H)$, then $F_\l:BH/H\to \C$, $b^{-1}H\mapsto \l(b)$ is a regular function on $BH/H$, and therefore a $B$-eigenvector of weight $\l$ in $\C(G/H)$.
In case $G/H$ is quasi affine (as for conjugacy classes), then $\Z\,\l(G/H)\cap P^+= \tilde\l(G/H)$ since $\C(G/H)=\Frac \C[G/H]$, as in \cite{pany4}, Proposition 1.5. I do not know if $\Z\,\l(G/H)\cap P^+= \tilde\l(G/H)$ holds in general.

\begin{lemma}\label{intero}
Suppose $F$ in $\Frac \C[G/H]$ is a $B$-eigenvector of weight $\l$ and $m\l$ lies in $\l(G/H)$ for a positive integer $m$. Then $F$ lies in $\C[G/H]$.
\end{lemma}
\pf There exists a $B$-eigenvector $F_1\in \C[G/H]$ of weight $m\l$. Then $F^m/F_1$ is
invariant under $B$ (as its weight is 0). So $F^m/F_1$ is constant,
as $G/H$ is spherical.  In other words, $F^m$ is regular on $G/H$. We conclude that $F$ is in $\C[G/H]$, since $\C[G/H]$ is integrally closed (\cite{grosshans}, Lemma 1.8).\cvd

\medskip

Let $\O$ be a spherical conjugacy class of $G$. We recall that $w=w(\O)$ is the unique element (an involution) of $W$ such that $BwB\cap \O$ is (open) dense in $\O$. Let $\v$ be the dense $B$-orbit in $\O$. Then $BG_y$ is dense in $G$ for any $y\in \v$. The parabolic subgroup $P=P_J$ associated to $G_y$ coincides with $\{g\in G\mid g.\v=\v\}$. Moreover $\v=\O\cap BwB$ (\cite{CCC}, Corollary 26), and it is affine, as an orbit of a soluble algebraic group.

We have $w=w_0\wJ$,  the subset $J$ is invariant under $\vartheta$, where $\vartheta$ is the symmetry of $\Pi$ induced by $-w_0$, and $w_0$ and $\wJ$ act in the same way on $\Phi_J$ 
(see \cite{Car} the discussion at the end of section 3, Corollary 4.2, Remark 4.3 and Proposition 4.15). 

Since all Levi subgroups of $P$ are conjugate under $R_uP$, we may choose $y\in\v$ such that the standard Levi subgroup $L_J$ is adapted to $G_y$. For the rest of this section
we fix such a $y$, and we put $H=G_y$, $P=P_J$, $L=L_J$. By Theorem \ref{nane2}, we have 
\begin{equation}
H_0=\<{U^-,U_{\wJ},C_y}=\<{U^-,U_{\wJ},T_y}
\label{nane4}
\end{equation}
and $\tilde\l(\O)=\l(G/H_0)$. 

We shall now relate $H$ with centralizers of elements in $\v\cap wB$.
 By the Bruhat decomposition, $y$ is of the form $y=u\dot w b$, where $u\in R_uP$ and $b\in B$. We put $x_1=u^{-1}yu=\dot wbu$. By \cite{Car}, Corollary 4.13, $U_{\wJ} (T^w)^\circ\leq C(x_1)$. Moreover, since $L'\leq C(y)$, by \cite{Car}, Lemma 3.4, and commutation of $y$ with $X_{\pm \a_i}$ for $i\in J$, we get $L'\leq C(x_1)$ (see also the proof of \cite{Car}, Proposition 4.15).

\begin{proposition}\label{nane5} Let $x$ be in $\O\cap wB$. Then
$
T_x=T_y
$ and $T\cap H^\circ=T\cap C(x)^\circ$.
\end{proposition}
\pf We observe that 
$C_{TU_w}(x)\leq T$ by the Bruhat decomposition and $C_{TU_w}(y)\leq T$, since $L$ is adapted to $C(y)$.
Now $x_1=u^{-1}yu=y^u$ implies
$$
\begin{array}{c}
T_{x_1}=C_T(x_1)=C_{TU_w}(x_1)\leq T\cap T^u = C_T(u)
\\
T_y=C_T(y)=C_{TU_w}(y)\leq T\cap T^{u^{-1}} = C_T(u^{-1})=C_T(u)
\end{array}
$$
therefore
if $t\in T_{y}$, then $t=t^u\in T_{x_1}$ and similarly if $t\in T_{x_1}$, then $t=t^{u^{-1}}\in T_y$. Hence
$
T_y=T_{x_1}
$, and $T\cap C(y)^\circ=T\cap C(x_1)^\circ$. To conclude note that  $\O\cap wB$ is the $T$-orbit of $x_1$.\cvd

\begin{remark}
{\rm 
In fact $C_L(x)=C_L(y)$ for every $x\in \O\cap wB$, since  $L'\leq C(x)$. 
}
\end{remark}

\begin{remark}
{\rm 
In general it is not true that $L_J$ is adapted to $C(x)$ for $x\in \O\cap wB$. For example if $\O$ is the minimal unipotent class, and $u$ is a non-identity element in $X_{-\beta}$, where $\beta$ is the highest root, then $C(u)\geq U^-$, so that there is a unique Levi subgroup of $P$ adapted to $C(u)$ (\cite{BP}, Proposition 3.9), and this is $L_J$. Since $u\not\in wB$, there is no element $x\in wB$ such that $L_J$ is adapted to $C(x)$. 
}
\end{remark}

From Theorem \ref{nane2} we get

\begin{corollary}\label{nane1} Let $\O$ be a spherical conjugacy class, $w=w(\O)$ and $x$ any element in $\O\cap wB$.
Then 
$
H_0=\<{U^-,U_{\wJ},T_x}
$, $w=w_0\wJ$.
\cvd
\end{corollary}

By Proposition \ref{nane5}, we may put $T_\O=T_x$, for $x\in \O\cap wB$. Then $T_\O=T_y$ and $(T^w)^\circ \leq T_\O\leq T^w$ by \cite{CCC}, step 2 in the proof of Theorem 5.

We shall need the description of the monoid of weights $\l$ such that $w(\l)=-\l$. In the next lemma we consider more generally $w$ of the form $w=w_0\wJ$, with $J$ $\vartheta$-invariant.

\begin{lemma}\label{descrizione} Let $J\subseteq \Pi$ be $\vartheta$-invariant and $w=w_0\wJ$. 
The dominant weight $\lambda$ satisfies $w(\lambda)=-\lambda$ if and only if $\lambda=\sum_{i\in \Pi\setminus J}n_i\omega_i$ with $n_{\vartheta(i)}=n_i$ for all $i\in \Pi\setminus J$. Moreover $w(\lambda)=-\lambda$ implies $w_0(\lambda)=-\lambda$.
\end{lemma}
\pf Let $\l\in P^+$, $\l=\sum n_i\omega_i$, $n_i\in\N$. For  $i\in \Pi\setminus J$ we have
$
\wJ(\omega_i)=\omega_i
$, so that $w(\omega_i)=-\omega_{\vartheta(i)}
$.

It is clear that if
$
\l=\sum_{i\in \Pi\setminus J}n_i\omega_i
$
with $n_i=n_{\vartheta(i)}$ for every $i\in \Pi\setminus J$, then
$
(w+1)(\l)=0$. 
On the other hand, assume $w(\lambda)=-\lambda$. Then $\wJ(\l)=-w_0\l$ and, by \cite{hum}, Theorem 1.12 (a), we get $-w_0\l=\l$ and $(\l,\a_j)=0$ for every $j\in J$. Hence $n_j=0$ for every $j\in J$. Moreover, from $\lambda=\sum_{i\in \Pi\setminus J}n_i\omega_i$ and
 $-w_0\l=\l$  it follows $n_{\vartheta(i)}=n_i$ for all $i\in \Pi\setminus J$.\cvd

\begin{remark}\label{base}{\rm 
If $S$ is a $\vartheta$-orbit in $\Pi\setminus J$, and we put
$
\omega_S=\sum_{i\in S}\omega_i
$
then we have seen that  $\{\omega_S\mid S\in (\Pi\setminus J)/\vartheta\}$ is a basis of the monoid $\{\l\in P^+\mid w(\l)=-\l\}$, where $(\Pi\setminus J)/\vartheta$ is the set of $\vartheta$-orbits in $\Pi\setminus J$. If we also assume that $w$ acts trivially on $\Phi_J$ (as in the case of $w=w(\O)$), then $\{\omega_S\mid S\in (\Pi\setminus J)/\vartheta\}$ is a basis of $\ker(w+1)$ in $E$, and so a basis of the free abelian group  $\{\l\in P\mid w(\l)=-\l\}$.
}\end{remark}

We describe $\tilde\l(\O)$. For this purpose we denote by $S_\O$ any supplement of $(T^w)^\circ$ in $T_\O$ (i.e. $S_\O(T^w)^\circ=T_\O$). We also put $P_w^+=\{\l\in P^+\mid w(\l)=-\l\}$. By Lemma \ref{descrizione} each element of $P_w^+$ satisfies $-w_0\l=\l$, so that in particular any subset $X$ of $P^+_w$ is {\em symmetric}, i.e. $-w_0(X)=X$ (\cite{pany3}, 4.2, \cite{Car}, Theorem 4.17)).

\begin{theorem}\label{main} Let $\O$ be a spherical conjugacy class, $w=w(\O)$ and let $S_\O$ be any supplement of $(T^w)^\circ$ in $T_\O$.
Then 
$$
\tilde\l(\O)=\{\l\in P^+_w\mid  \l(S_\O)=1\}
$$
\end{theorem}
\pf By (\ref{questo}), $\tilde\l(\O)=\{\l\in P^+\mid \l(T_\O)=1\}$. Since $(T^w)^\circ\leq T_\O$,  a necessary condition for $\l\in P^+$ to be in $\tilde\l(\O)$ is that  $\l(t\,t^w)=1$ for every $t\in T$, as $(T^w)^\circ=\{t\,t^w\mid t\in T\}$. This condition is equivalent to $(w+1)\l=0$, so that $\tilde\l(\O)\leq P^+_w$. Let $\l\in P^+_w$: then  $\l\in\tilde\l(\O)\Iff \l(S_\O)=1$.
\cvd

\medskip

We shall prove the crucial fact that $\tilde\l(\O)=\l(\O)$, so that the monoid $\l(\O)$ is {\em saturated} (that is $\Z\,\l(\O)\cap P^+=\l(\O)$, \cite{pany4}, Definition 1.3). In the following, $x$ is a fixed element in $\O\cap wB$ and $\w$ a representative of $w$ in $N$ such that $x=\w u$, $u\in U$. If $u=\prod_{\a\in\Phi^+}x_\a(k_\a)$, and $i\in \Pi$, we say that $\a_i$ {\em occurs} in $x$ if $k_{\a_i}\not=0$. This is independent of the chosen total ordering on $\Phi^+$.

For the closure $\ov\O$ of $\O$ in $G$, the monoid $\l(\ov\O)$ of dominant weights occurring in $\C[\ov\O]$ is a submonoid of $\l(\O)$. We start with

\begin{proposition}\label{funzioni} Let $\l\in P^+$. Then $(1-w)\l$ lies in $\l(\overline\O)$.
\end{proposition}
\pf Let $f\in V(\l)^\ast_{-w\l}$, $v\in V(\l)_\l$ with $f(\dot w. v)=1$. Then
$
c_{f,v}(t^{-1}gt)= c_{t.f,t.v}(g)=((1-w)\l)(t)c_{f,v}(g)
$
for every $t\in T$, $g\in G$.
For every $z$, $z_1\in U$ we have
$$
c_{f,v}(z_1xz)= f(z_1\dot w\, u z.v)=f(z_1\dot w\,.v)=f(\dot w\,.v)=1
$$
since $z_1\dot w\,.v=\dot w\,.v+v_1$, where $v_1$ is a sum of weight vectors of weights strictly greater than $w\l$. Therefore for every $t\in T$, $z\in U$ we have
\begin{equation}
c_{f,v}(t^{-1}z^{-1}xzt)=((1-w)\l)(t)
\label{peso}
\end{equation}
Since $B.x$ is dense in $\ov\O$, by (\ref{peso})
the restriction of $c_{f,v}$ to $\ov\O$ is a (non-zero) $B$-eigenvector of weight $(1-w)\l$ in $\C[\ov\O]$. Hence $(1-w)\l\in \l(\ov\O)$.\cvd

\begin{corollary}\label{funzioni3} Let $\l\in P^+_w$. Then $2\l$ lies in  $\l(\overline\O)$.
\cvd
\end{corollary}

\begin{corollary}\label{funzioni4} Let $\l\in P^+$. Then $(1-w)\l\in\l(\O)$. If moreover $\l\in P^+_w$, then $2\l$ lies in  $\l(\O)$.
\end{corollary}
\pf This follows from the fact that $\l(\overline\O)\leq \l(\O)$.\cvd

\medskip
We have shown that
\begin{equation}
2P^+_w\leq (1-w)P^+\leq\l(\ov\O)\leq \l(\O)\leq \tilde\l(\O)\leq P_w^+\label{catena}
\end{equation}
We can prove that $\l(\O)$ is saturated.

\begin{theorem}\label{saturo} Let $\O$ be a spherical conjugacy class. Then $\l(\O)$ is saturated.
\end{theorem}
\pf 
 Let $\l\in \tilde\l(\O)$. We put
$
F(b^{-1}xb)=\lambda(b)
$
for $b\in B$. We observed that  $F$ is well-defined since $C_B(x)=T_xU_{\wJ}$ and gives rise to a $B$-eigenvector of weight $\l$ in $\C(\O)$. Since $\O$ is quasi affine, we conclude that $\l$ lies in $\l(\O)$ by Theorem \ref{main}, Corollary  \ref{funzioni4} and Lemma \ref{intero}.\cvd

Theorem \ref{saturo} in particular proves Conjecture 5.12 (and 5.10 and 5.11) in \cite{pany5}.

\medskip

To deal with $\l(\ov\O)$, in section \ref{non-normal} we shall make use of

\begin{proposition}\label{funzioni5} Let $\l\in P^+$, $i\in \Pi\setminus J$ be such that $\a_i$ occurs in $x$ and $(\l,\a_i)\not=0$. Then $(1-w)\l-\a_i\in  \l(\overline\O)$.
\end{proposition}
\pf Since $\<{\l,\a_i}\not=0$, $\l-\a_i$ is a weight of $V(\l)$. We construct two matrix coefficients.  We fix a non-zero $v\in V(\l)_{\l-\a_i}$. By \cite {yale},  Lemma 72, there exists a (unique) $v_\l\in V(\l)_{\l}$ such that 
$
x_{\a_i}(k).v=v+kv_\l
$
for every $k\in \C$. Then we choose $f\in V(\l)^\ast_{-w\l}$ such that 
$f(\dot w.v_\l)=1$.

Since $\a_i$ occurs in $x=\w\,u$, we have $u=x_{\a_i}(r)u'$, with $r\in \C^\ast$, $u'\in \prod_{\b\in\Phi^+\setminus \{\a_i\}}X_\b$. Let  $y$, $y_1\in U$, and let  $y=x_{\a_i}(k)y'$, $y'\in  \prod_{\b\in\Phi^+\setminus \{\a_i\}}X_\b$, then
$$
y_1^{-1}xy.v=
y_1^{-1}\dot w.v+(k+r) y_1^{-1}\dot w.v_\l
$$
The vector $\dot w.v$ has weight $w(\l-\a_i)$, so that $y_1^{-1}\dot w.v$ is a sum of weight vectors of weight $w(\l-\a_i)+\beta$, where $\beta$ is a sum of simple roots with non-negative coefficients. Assume $w\l=w(\l-\a_i)+\beta$ for a certain $\beta$. Then $w(\a_i)=\beta$ would be positive, a contradiction since $i\in \Pi\setminus J$. Hence $f(y_1^{-1}\dot w.v)=0$. Similarly, 
$y_1^{-1}\dot w.v_\l=\dot w.v_\l+v'$, where $v'$  is a sum of weight vectors of weights greater than $w\l$, hence $f(y_1^{-1}\dot w.v_\l)=f(\dot w.v_\l)=1$, so that 
$
c_{f,v}(y_1^{-1}xy)= k+r
$.

The second matrix coefficient is defined dually. We fix a non-zero $f_1\in V(-w_0\l)^\ast_{\l-\a_i}$.  There exists a (unique) $f_\l\in V(-w_0\l)^\ast_{\l}$ such that 
$
x_{\a_i}(k).f_1=f_1+kf_\l
$
for every $k\in \C$. Then we choose $v_1\in V(-w_0\l)_{-w\l}$ such that 
$f_\l(\dot w.v_1)=1$. 
Let $z$, $z_1\in U$,   $z_1=x_{\a_i}(k_1)z'$, $z'\in  \prod_{\b\in\Phi^+\setminus \{\a_i\}}X_\b$, then proceeding as before, we get
$
c_{f_1,v_1}(z_1^{-1}xz)= k_1
$.

For $t\in T$, $z\in U$ we obtain
\begin{equation}
(c_{f,v}-c_{f_1,v_1})(t^{-1}z^{-1}xzt)= r\,((1-w)\l-\a_i)(t)
\label{riassunto}
\end{equation}
Since $B.x$ is dense in $\ov\O$, by (\ref{riassunto})
the restriction of $c_{f,v}-c_{f_1,v_1}$ to $\ov\O$ is a (non-zero) $B$-eigenvector of weight $(1-w)\l-\a_i$ in $\C[\ov\O]$. Hence $(1-w)\l-\a_i\in \l(\ov\O)$.\cvd

\begin{corollary}\label{funzioni18} Let $i\in \Pi\setminus J$ be such that $\a_i$ occurs in $x$. Then $\omega_i+\omega_{\vartheta(i)}-\a_i$ lies in $\l(\overline\O)$.
\end{corollary}
\pf This follows from Proposition \ref{funzioni5} by taking $\l=\omega_i$.\cvd

\medskip
We can deal with other homogeneuos spaces related to $\O$. The simply-connected cover (or the universal covering, as in \cite{Jan}, p. 107) $\hat\O$ of $\O$ can be identified with $G/H^\circ$, since $G$ is simply-connected.

\begin{corollary}\label{covering} Let $\O$ be a spherical conjugacy class, and let $S$ be a supplement of $(T^w)^\circ$ in $T\cap C(x)^\circ$. Then $\l(\hat\O)=\{\l\in P^+_w\mid \l(S)=1\}$ is saturated.
\end{corollary}
\pf By \cite{grosshans}, Corollary 2.2, $\hat\O$ is quasi affine and, by \cite{BP}, Proposition 5.1, 5.2, $L$ is adapted to $H^\circ$, so that $\tilde\l(\hat\O)=\tilde\l(G/H^\circ)=\{\l\in P^+_w\mid \l(S)=1\}$, since $(T^w)^\circ\leq T\cap H^\circ$. Let  
$\l\in\tilde\l(\hat\O)$; 
then $F_\l:BH^\circ/H^\circ\to \C$, $b^{-1}H^\circ\mapsto \l(b)$ is a regular function on $BH^\circ/H^\circ$, and therefore a $B$-eigenvector of weight $\l$ in $\C(G/H^\circ)$.
By Corollary \ref{funzioni4}, $2\l\in \l(G/H)\leq \l(G/H^\circ)$, and we conclude by Lemma \ref{intero} and Proposition \ref{nane5}.\cvd

\begin{corollary}\label{Gen. saturo} 
Let $K$ be a closed subgroup of $G$ with $H^\circ\leq K\leq N(H^\circ)$. Then $\l(G/K)=\tilde\l(G/K)$ (and
$\l(G/K)$ is saturated).
\end{corollary}
\pf Since $L$ is adapted to $H$, we get $N(H)=N(H^\circ)=H(C\cap N(H))$ by \cite{BP}, Corollaire 5.2, $P$ is the parabolic subgroup corresponding to $N(H)$ and $L$ is adapted to $N(H)$ (by the proof of \cite{BP}, Proposition 5.2 a). Clearly the same holds for $K$, since $BH=BK$.

By Corollary \ref{covering},  $\l\in \l(G/H^\circ)\iff \l(T\cap H^\circ)=1$. We prove that $\l\in \l(G/K)\iff \l(T\cap K)=1$. In one direction $\l\in \l(G/K)\Ra \l(T\cap K)=1$, since $\l(G/K)\leq \tilde \l(G/K)$. So assume $\l(T\cap K)=1$.
Then $\l(T\cap H^\circ)=1$, so that $\l\in \l(G/H^\circ)$, and in particular $w_0\l=-\l$. Let $v$ be a non-zero vector in $V(\l)^{H^0}$, and let $v=v_{-\l}+v'$, with $v_{-\l}\in V(\l)_{-\l}$, $v'\in \sum_{\mu>-\lambda} V(\l)_\mu$: then $v_{-\l}\not=0$, since $BH^\circ$ is dense in $G$. 

Since $V(\l)^{H^0}$ is 1-dimensional, there  is a character $\gamma$ of $K$, trivial on $H^\circ$, such that  $k.v=\gamma(k)v$ for $k\in K$.
 Since $K=H^\circ(T\cap K)$,  $v$ is $K$-invariant if and only if $\gamma(T\cap K)=1$. But $v_{-\l}\not=0$ implies
$\gamma(k)=-\l(k)$ for every $k\in T\cap K$ so that $v$ is $K$-invariant if and only if $\l(T\cap K)=1$, and we are done.
\cvd

\begin{remark} 
{\rm
In general $K$ is not quasi affine: for instance the centralizer $H$ of $x_{-\b}(1)$, $\b$ the highest root, contains $U^-$, and $T\leq N(H)$. Then $N(H)$ is epimorphic, i.e. the minimal quasi affine subgroup of $G$ containing $N(H)$ is $G$ (\cite{grosshans}, p. 19, ex. 2).  To our knowledge, it was known that $\l(G/K)$ is saturated for symmetric varieties $G/K$, due to the work of Vust, \cite{Vust}.}
\end{remark}

\begin{proposition}\label{T e G} We have
$$
H/H^\circ\cong T_y/T\cap H^\circ  =
 T_x/T\cap C(x)^\circ
$$
\end{proposition}
\pf We have $H=H^\circ(H\cap T)=H^\circ T_y$. Hence we get an epimorphism $\pi:T_y\to H/H^\circ$, inducing an isomorphism $\ov\pi: T_y/T\cap H^\circ\to H/H^\circ$, and we conclude by Proposition \ref{nane5}.\cvd

\begin{corollary}\label{connesso} If $T^w$ is connected, then $H$ is connected.
\end{corollary}
\pf This follows from $(T^w)^\circ\leq T\cap C(x)^\circ\leq T_x\leq T^w=(T^w)^\circ$ and Proposition \ref{T e G}.\cvd
\medskip

Due to the fact that $T$ is $2$-divisible, we have the decomposition
$
T=(T^w)^\circ (S^w)^\circ
$
where $S^w=\{t\in T\mid t^w=t^{-1}\}$. Let $t\in T^w$,  $t=s\,z$, with $s\in (T^w)^\circ$, $z\in (S^w)^\circ$. Then $z=t\,s^{-1}\in T^w\cap (S^w)^\circ\leq T^w\cap S^w\leq T_2$, 
the elementary abelian $2$-subgroup of $T$ of rank $n$. We note that $(T^w)^\circ\cap (S^w)^\circ$ is finite, even though 
in general  not trivial. Therefore $z\in T_2$, and $T^w\leq (T^w)^\circ\, T_2$. In particular we have 
$$
T^w= (T^w)^\circ (T^w\cap (S^w)^\circ)=(T^w)^\circ (T^w\cap T_2)
$$
and
$$
T_x= (T^w)^\circ (C(x)\cap (S^w)^\circ)=(T^w)^\circ (C(x)\cap T_2)
$$
Moreover every subgroup $M$ of $T_2$ is a complemented group (i.e. for every subgroup $X$ of $M$ there exists a subgroup $Y$ such that $X\,Y=M$ and $X\cap Y=1$), hence we may find a subgroup $R$ of $T_2$ such that $T^w=(T^w)^\circ\times R$. Then   $T_x=(T^w)^\circ\times (R\cap C(x))$ and
$T\cap C(x)^\circ=(T^w)^\circ\times (R\cap C(x)^\circ)$. We put
$S_\O=R\cap C(x)$, $S_{\hat\O}=R\cap C(x)^\circ$. We have therefore proved

\begin{theorem}\label{main-bis} Let $\O$ be a spherical conjugacy class, $w=w(\O)$.
Then 
$$
\l(\O)=\{\l\in P^+_w\mid  \l(S_\O)=1\}
\quad ,\quad 
\l(\hat\O)=\{\l\in P^+_w\mid  \l(S_{\hat\O})=1\}
$$
\cvd
\end{theorem}

From Proposition \ref{T e G} it follows that $H$ always splits over $H^\circ$: if $Y$ is a complement of $R\cap C(x)^\circ$ in $R\cap C(x)$, then $Y$ is a complement of $H^\circ$ in $H$.

\section{Description of $\l(\O)$ and $\l(\hat \O)$}\label{vari}

From our discussion it is clear that to determine $\l(\O)$ the most favourable case is when $T^w$ is connected, so that $T_x=T^w=(T^w)^\circ$. In this case then $\l(\O)=\l(\hat\O)=P^+_w=\{\sum_{i\in\Pi\setminus J}n_i\omega_i\mid n_{\vartheta(i)}=n_i\}$.
We note that of course we have $Z(G)\leq T_x$, so that it is also straightforward to determine $\l(\O)$ even when $T^w=(T^w)^\circ Z(G)$, so that $T_x=T^w$. In general it is quite cumbersome to determine $T_x$. Our strategy will be to determine $T^w$ as $T^w=(T^w)^\circ \times R$, and then determine $R\cap C(x)$. To deal with unipotent classes, we shall usually start from the maximal one, (corresponding to $w_0$), and then deal with the remaining classes by an inductive procedure. In some cases we shall use an explicit form of an element $x$ (in $\O\cap wB$), while in some other cases we shall determine $T\cap C(x)$ by analizing the form of eventual involutions in $T_x\setminus Z(G)(T^w)^\circ$.  Note that when $T^w$ is connected (or $T^w=(T^w)^\circ Z(G)$), it is not necessary to have an explicit description of $x\in \O\cap wB$
(however in certain cases it will be necessary to have such a description in section \ref{generale}).

We use the fact that if $G_1\subset G_2$ are reductive
algebraic groups and $u$ is a unipotent element in $G_1$ such that the conjugacy class of  $u$ in $G_2$ is spherical, then the
conjugacy class of $u$ in $G_1$ is spherical (\cite{pany}, Corollary 2.3,
Theorem 3.1). 

The character group $X(T^w)$ is isomorphic to $P/(1-w)P$, since $P=X(T)$. Therefore $T^w$ is connected if and only if $P/(1-w)P$ is torsion free. We are reduced to calculate elementary divisors of the endomorphism $1-w$ of $P$.  We shall use the following results.

\begin{lemma}\label{scambio} Assume the positive roots $\b_i,\ldots,\b_\ell$ are long and pairwise orthogonal. Then, for $\xi_1,\ldots,\xi_\ell\in \C^\ast$ and  $g=x_{\b_1}(-\xi_1^{-1})\cdots x_{\b_\ell}(-\xi_\ell^{-1})$ we have
$$
gx_{-\b_1}(\xi_1)\cdots x_{-\b_\ell}(\xi_\ell)g^{-1}=
n_{\b_1}\cdots n_{\b_\ell}hx_{\b_1}(2\xi_1^{-1})\cdots x_{\b_\ell}(2\xi_\ell^{-1})
$$
for a certain $h\in T$.
\end{lemma}
\Pf By (\ref{relazioni}) we have
$
x_\a(-\xi^{-1})x_{-\a}(\xi)x_\a(\xi^{-1})=n_\a h_\a(-\xi)x_\a(2\xi^{-1})
$. 
Hence we get the result with
 $h=h_{\b_1}(-\xi_1)\cdots h_{\b_\ell}(-\xi_\ell)$.
\cvd

\begin{proposition}\label{minima} Let $\a\in\Phi$. Then $T^{s_\a}$ is connected except in the following cases: 
\begin{itemize}
\item[(i)]  $G$ is of type $A_1$;
\item[(ii)] $G$ is of type $C_n$ and $\a$ is long;
\item[(iii)] $G$ is of type $B_2$ and $\a$ is long.
\end{itemize}
In these cases we have $T^{s_\a}= (T^{s_\a})^\circ\times Z(G)$.
\end{proposition}
\pf 
It is enough to determine in which cases the non-zero elementary divisor of $1-s_i$ is not 1.  Since $(1-s_i)\omega_j=\d_{ij}\a_i$ and $\a_i=\sum_k a_{ik}\omega_k$, this happens only for $G$ of type $A_1$ and $i=1$, $C_n$ and $i=n$, or $B_2$ and $i=1$  (\cite{Hu}, pag. 59). In these cases the non-zero elementary divisor is 2, and $T^{s_{\a_i}}= (T^{s_{\a_i}})^\circ\times Z(G)$. \cvd

\begin{lemma}\label{esclude} Let $M$ be a connected algebraic group, $S$ a torus of $M$, $g$ a semisimple element in $C_M(S)$. Then $\<{S,g}$ is contained in a torus of $M$.
\end{lemma}
\pf See \cite{Hu}, Corollary 22.3 B. \cvd

\begin{lemma}\label{free}
Assume $K$ is a connected spherical subgroup of $G$ with no non-trivial characters. Then the monoid $\l(G/K)$ is free.
\end{lemma}
\pf We recall that we are assuming $G$ simply-connected, so that by \cite{grosshans}, Theorem 20.2, ${}^U\C[G/K]$ is a polynomial algebra. But ${}^U\C[G/K]$ is the monoid algebra of $\l(G/K)$ and the monoid algebra is factorial if and only if $\l(G/K)$ is free (see the proof of \cite{pany3}, Proposition 2).
\cvd

\begin{lemma}\label{dim 1}
Let $V$ be a $G$-module, $g\in G$, such that the image $Q$ of the endomorphism $p(g)$ of $V$ is 1 dimensional for a certain polynomial $p$.
Assume $M\leq C(g)$ has no non-trivial characters. Then $M$ acts trivially on $Q$.
 \end{lemma}
\pf This is clear.\cvd 

\medskip
\noindent
Let $S=\{i,\vartheta(i)\}$ be a $\vartheta$-orbit in $\Pi\setminus J$ consisting of 2 elements. We put $H_S=\{h_{\a_i}(z)h_{\a_{\vartheta(i)}}(z^{-1})\mid z\in \C^\ast\}
$. Let ${\cal S}_1$ be the set of $\vartheta$-orbits in $\Pi\setminus J$ consisting of 2 elements.  Then, by Remark \ref{base}, $\Delta_J\cup \{\a_i-\a_{\vartheta(i)}\}_{{\cal S}_1}$ is a basis of $\ker(1-w)$ and
\begin{equation}
(T^w)^\circ=\prod_{j\in J}H_{\a_j}\times \prod_{S\in {\cal S}_1}
H_S
\label{prodotto}
\end{equation}
We put $\Psi_J=\{\beta\in \Phi\mid w(\b)=-\b\}$. Then $\Psi_J$ is a root system in $\im(1-w)$  (\cite{springer2}, Proposition 2), and $w_{|\im(1-w)}$ is $-1$.
If $K=C((T^w)^\circ)'$, then $K$ is semisimple with root system $\Psi_J$ and maximal torus $T(K):=T\cap K=(S^w)^\circ$.

For each spherical  (non-central) conjugacy class $\O$ we give the corresponding $J$ and $w$ as a product of commuting reflections using the tables in \cite{CCC}.  We give tables with corresponding $\l(\O)$ and $\l(\hat\O)$ (for semisimple classes we also give the type of the centralizer of elements in $\O$). In the cases when $\l(\hat\O)=\l(\O)$, we leave a blank entry. For length reasons we shall give proofs only for some classes. In \cite{CCC} for the classical groups we gave representative of semisimple conjugacy classes in $SL(n)$, $Sp(n)$ and $SO(n)$. Here we shall give an expression in terms of $\exp$. If $g$ is in $Z(G)$, then $\O_g=\{g\}$, $w=1$ and $\C[\O_g]=\C$.

\medskip

\subsection{Type $A_{n}$, $n\geq 1$.}

Let $m=\left[\frac{n+1}{2}\right]$, $\b_i=e_{i}-e_{n+2-i}$, for $i=1,\ldots, m$. 
For $\ell=1,\ldots,m-1$ we put $J_\ell=\{\ell+1,\ldots,n-\ell\}$, $J_m=\vuoto$. 

\medskip

\subsubsection{Unipotent classes in $A_n$.} 

If we denote by $X_i$ the unipotent class $(2^i,1^{n+1-2i})$, then 
$$
X_\ell \corr J_\ell \corr s_{\b_1}\cdots s_{\b_\ell}
$$
for $\ell=1,\ldots,m$ (here $w_0=s_{\b_1}\cdots s_{\b_m}$). 

In this case $T^w$ is almost always connected. There is only one case when it is not connected, namely when $n$ is odd, $n+1=2m$, and $w=w_0$. However in this case we have $T^{w_0}=(T^{w_0})^\circ Z(G)=(T^{w_0})^\circ \times\<{h_{\a_m}(-1)}$.

We get
\begin{center}
\vskip-20pt
$$
\begin{array}{|c||c|c|}
\hline
\O  & \l(\O)& \l(\hat\O) \\ 
\hline
\hline
\begin{array}{c}
X_\ell\\
\ell=1,\ldots,m-1
\end{array}
& \displaystyle\sum_{k=1}^\ell n_k (\omega_k+\omega_{n-k+1})&  \\
\hline
\begin{array}{c}
X_m\\
n=2m
\end{array}
& \displaystyle\sum_{k=1}^m n_k (\omega_k+\omega_{n-k+1})&  
\\
\hline
\begin{array}{c}
X_m\\
n+1=2m
\end{array}
&  \displaystyle\sum_{k=1}^{m-1} n_k (\omega_k+\omega_{n-k+1})+2n_m\omega_m&  \displaystyle\sum_{k=1}^{m-1} n_k (\omega_k+\omega_{n-k+1})+n_m\omega_m\\
\hline
\end{array}
$$
\end{center}
\begin{center} Table \totable: $\l(\O)$, $\l(\hat\O)$ for unipotent classes in $A_n$.
\end{center}
In particular $\hat X_1$ is a model homogeneus space for $SL(2)$, and in fact the principal one, by \cite{Luna}, 3.3 (1).

\medskip
\subsubsection{Semisimple classes in $A_n$.} 

Following the notation in \cite{CCC}, Tables 1, 5 we get
$$
T_1A_{\ell-1}A_{n-\ell} \corr J_\ell \corr s_{\b_1}\cdots s_{\b_\ell}
$$
for $\ell=1,\ldots,m$.

We get

\begin{center}
\vskip-20pt
$$
\begin{array}{|c|c||c|}
\hline
\O  & H & \l(\O) \\ 
\hline
\hline
\begin{array}{c}
\exp(\zeta\check\omega_\ell)\\
\zeta\in \C\setminus 2\pi i\Z\\
\ell=1,\ldots,m-1
\end{array}
& T_1 A_{\ell-1}A_{n-\ell}&  \displaystyle\sum_{k=1}^\ell n_k (\omega_k+\omega_{n-k+1})\\
\hline
\begin{array}{c}
\exp(\zeta\check\omega_m)\\
\zeta\in \C\setminus 2\pi i\Z\\
n=2m
\end{array}
& T_1 A_{m-1}A_{m}&  \displaystyle\sum_{k=1}^m n_k (\omega_k+\omega_{n-k+1})\\
\hline
\begin{array}{c}
\exp(\zeta\check\omega_m)\\
\zeta\in \C\setminus 2\pi i\Z\\
n+1=2m
\end{array}
& T_1 A_{m-1}A_{m-1}&  \displaystyle\sum_{k=1}^{m-1} n_k (\omega_k+\omega_{n-k+1})+2n_m\omega_m\\
\hline
\end{array}
$$
\end{center}
\begin{center} Table \totable: $\l(\O)$ for semisimple classes in $A_n$.
\end{center}

\subsection{Type $C_n$, $n\geq 2$.}

We have $\omega_\ell=e_1+\cdots+e_\ell$ for $\ell=1,\ldots,n$ and
$
Z(G)=\<z$,
where 
$
z=\prod_{i=1}^{[\frac {n+1}2]}h_{\a_{2i-1}}(-1)
$.

\subsubsection{Unipotent classes in $C_n$.} 

For $i=1,\ldots,n$
we denote by $X_i$ the unipotent class $(2^{i},1^{2n-2i})$ and we put
$\b_i=2e_{i}$, $J_i=\{i+1,\ldots,n\}$ ($J_n=\vuoto$).
Then
$$
X_\ell \corr J_\ell \corr s_{\b_1}\cdots s_{\b_\ell}
$$
for $\ell=1,\ldots,n$ (here $w_0=s_{\b_1}\cdots s_{\b_n}$).

\begin{lemma}\label{toni1c} Let $w=s_{\b_1}\cdots s_{\b_\ell}
$
for $\ell=1,\ldots,n$. Then 
$$
T^w=(T^w)^\circ \times R\quad,\quad
R=\<{h_{\a_1}(-1)}\times\cdots\times \<{h_{\a_\ell}(-1)}
$$
\end{lemma}
\pf For $\ell=1,\ldots,n$ we have $(1-w)P=\Z\<{2\omega_1,\ldots,2\omega_\ell}$.\cvd

\begin{proposition}\label{toni2c} For $\ell=1,\ldots,n$ we have
$$
\l(X_\ell)=\{2n_1\omega_1+\cdots+2n_\ell\omega_\ell\mid n_k\in \N\}
$$
\end{proposition}
\pf
In \cite{CCC} we exhibit the element  $x_{-\b_1}(1)\cdots x_{-\b_\ell}(1)\in \O\cap BwB\cap B^-$. By Lemma \ref{scambio}, we can choose
$$
x=n_{\b_1}\cdots n_{\b_\ell}h\,x_{\b_1}(2)\cdots x_{\b_\ell}(2)
\in \O\cap wB
$$
for a certain $h\in T$. Let now $t\in R$. Then $t\in C(x)\iff \b_i(t)=1$ for $i=1,\ldots,\ell$. But $\Z\<{\b_1,\ldots,\b_\ell}=\Z\<{2\omega_1,\ldots,2\omega_\ell}$, so that $R\leq T_x$, and $T_x=T^w$.\cvd

\begin{proposition}\label{toni2ccovering} For $\ell=1,\ldots,n$ we have
$$
\l(\hat X_\ell)=\{2n_1\omega_1+\cdots+2n_{\ell-1}\omega_{\ell-1}+n_\ell\omega_\ell\mid n_k\in \N\}
$$
\end{proposition}
\pf We have $R\cap C(x)^\circ=\<{h_{\a_1}(-1),\ldots,h_{\a_{\ell-1}}(-1)}$.
In fact, for $i=1\ldots, \ell-1$
$$
e_{\a_i}-e_{-\a_i}\in C_{\frak g}(\<{x_{\b_1}(\xi)\cdots x_{\b_\ell}(\xi)})
$$
for every $\xi\in \C$,
so that $h_{\a_i}(-1)=\exp(\pi (e_{\a_i}-e_{-\a_i}))\in C(x)^\circ$. On the other hand the reductive part of $C(x)$ is of type $Sp(2n-2\ell)\times O(\ell)$, so that $C(x)/C(x)^\circ$ has order 2, and we are done.\cvd

Hence

\begin{center}
\vskip-20pt
$$
\begin{array}{|c||c|c|}
\hline
\O  & \l(\O)& \l(\hat\O) \\ 
\hline
\hline
\begin{array}{c}
X_\ell\\
\ell=1,\ldots,n
\end{array}
&\quad \displaystyle\sum_{i=1}^\ell 2n_i\omega_i \quad&\quad\displaystyle\sum_{i=1}^{\ell-1} 2n_i\omega_i+n_\ell\omega_\ell\quad \\
\hline
\end{array}
$$
\end{center}
\begin{center} Table \totable: $\l(\O)$, $\l(\hat\O)$ for unipotent classes in $C_n$.
\end{center}

\medskip

\subsubsection{Semisimple classes in $C_n$.} 
Let $p=[\frac n2]$. We put  $\g_\ell=e_{2\ell-1}+e_{2\ell}$, $K_\ell=\{1,3,\ldots,2\ell-1,2\ell+1,2\ell+2,\ldots,n\}$
for $\ell=1,\ldots, p$. 
Then, following the notation in \cite{CCC}, Tables 1, 5 we have
$$
\begin{array}{lllll}
C_\ell C_{n-\ell},\quad \ell=1,\ldots,p&\corr &K_\ell&\corr &s_{\gamma_1}\cdots s_{\gamma_\ell}\\
T_1 C_{n-1}&\corr &J_2 &\corr&s_{\b_1}s_{\b_2}\\
T_1\tilde A_{n-1}&\corr &\vuoto &\corr & w_0
\end{array}
$$

We get

\begin{center}
\vskip-20pt
$$
\begin{array}{|c|c||c|}
\hline
\O  & H & \l(\O) \\ 
\hline
\hline
\begin{array}{c}
\exp(\zeta\check\omega_n)\\
\zeta\in\C\setminus 2\pi i\Z
\end{array}
& T_1 \tilde A_{n-1}&  \displaystyle\sum_{k=1}^n 2n_k \omega_k\\
\hline
\begin{array}{c}
\exp(\zeta\check\omega_1)\\
\zeta\in\C\setminus \pi i\Z
\end{array}
& T_1 C_{n-1}&  2n_1 \omega_1+n_2\omega_2\\
\hline
\begin{array}{c}
\exp(\pi i \check\omega_\ell)\\
\ell=1,\ldots,[{\frac n2}]
\end{array}
&C_\ell C_{n-\ell}&  \displaystyle\sum_{i=1}^\ell n_{2i}\,\omega_{2i}\\
\hline
\end{array}
$$
\end{center}
\begin{center} Table \totable: $\l(\O)$ for semisimple classes in $C_n$.
\end{center}

\medskip
\subsubsection{Mixed classes in $C_n$.}
We put $p=[\frac n2]$. From \cite{CCC}, Table 4, we get
$$
\begin{array}{llllll}
\s_p x_{\a_n}(1)&& \corr & \vuoto& \corr &  w_0\\
\s_k x_{\a_n}(1),\ &k=1,\ldots,p-1& \corr & J_{2k+1}& \corr &  s_{\b_1}\cdots s_{\b_{2k+1}}\\
\s_k x_{\b_1}(1),\ &k=1,\ldots,p& \corr & J_{2k}& \corr & s_{\b_1}\cdots s_{\b_{2k}}\end{array}
$$
\noindent
Note that when $n$ is even, then $\s_p x_{\b_1}(1)\sim z\s_p x_{\a_n}(1)$.
 We obtain

\begin{center}
\vskip-20pt
$$
\begin{array}{|c||c|c|}
\hline
\O  & \l(\O)& \l(\hat\O) \\ 
\hline
\hline 
\begin{array}{c}
\vspace{0.06cm} 
\s_p x_{\a_n}(1) 
\end{array}
& \displaystyle\sum_{i=1}^n n_i\omega_i,
\  \text{\ $\sum_{i=1}^{[\frac {n+1}2]}n_{2i-1}$ even} & \displaystyle\sum_{i=1}^n n_i\omega_i \\
\hline
\begin{array}{c}
\s_\ell x_{\a_n}(1)\\
\ell=1,\ldots,[\frac n2]-1
\end{array}
&\quad \displaystyle\sum_{i=1}^{2\ell+1} n_i\omega_i,
\  \text{\ $\sum_{i=1}^{\ell+1}n_{2i-1}$ even}\quad&\quad \displaystyle\sum_{i=1}^{2\ell+1} n_i\omega_i\quad
 \\
\hline
\begin{array}{c}
\s_\ell x_{\b_1}(1)\\
\ell=1,\ldots,[\frac n2]
\end{array}
&  \displaystyle \sum_{i=1}^{2\ell} n_i\omega_i
,\  \text{\ $\sum_{i=1}^{\ell}n_{2i-1}$ even}& \displaystyle \sum_{i=1}^{2\ell} n_i\omega_i
\\
\hline
\end{array}
$$
\end{center}
\begin{center} Table \totable: $\l(\O)$, $\l(\hat\O)$ for mixed classes in $C_n$.
\end{center}
In particular $\hat\O_{\s_p x_{\a_n}(1)}$ is a model homogeneus space, and in fact the principal one, by \cite{Luna}, 3.3 (3).

\medskip

To deal with types $D_n$ and $B_n$, we denote by $X_i$ the unipotent class which in $SO(s)$ has canonical form $(2^{2i},1^{s-4i})$, $i=1,\ldots,\left[\frac{s}{4}\right]$ (for $s=4m$, $i=m$ there are 2 classes of this form: $X_m$ and $X_m'$, the {\em very even} classes) and by $Z_i$ the unipotent class $(3,2^{2(i-1)},1^{s-4i+1})$, $i=1,\ldots,1+\left[\frac{s-3}{4}\right]$.

\subsection{Type $D_n$, $n\geq 4$.}
Let $m=\left[\frac{n}2\right]$. We have $\omega_i=e_1+\cdots+e_i$ for $i=1,\ldots, n-2$, $\omega_{n-1}=\frac12(e_1+\cdots+e_{n-1})-\frac 12 e_n$, $\omega_n=\frac12(e_1+\cdots+e_{n})$. 
We put $\b_i=e_{2i-1}+e_{2i}$, $\d_i=e_{2i-1}-e_{2i}$ for $i=1,\ldots, m$.
For $\ell=1,\ldots,m-1$ we put $J_\ell=\{2\ell+1,\ldots,n\}$, $J_m=\vuoto$, 
$K_\ell=J_\ell\cup\{1,3,\ldots,2\ell-1\}$ for $\ell=1,\ldots,m$.

\subsubsection{Unipotent classes in $D_n$, $n$ {even}, $n=2m$.} 

The center of $G$ is
$
\<{\prod_{i=1}^{m}h_{\a_{2i-1}}(-1),h_{\a_{n-1}}(-1)h_{\a_n}(-1)}$.
From \cite{CCC} we get
$$
\begin{array}{llllll}
Z_\ell,& \ell=1,\ldots,m&\corr &J_\ell&\corr &s_{\b_1}s_{\d_1}\cdots s_{\b_\ell}s_{\d_\ell}\\
X_\ell,& \ell=1,\ldots,m&\corr &K_\ell&\corr & s_{\b_1}\cdots s_{\b_\ell}\\
X_m' &&\corr & \{1,3,\ldots,n-3,n\}&\corr & s_{\b_1}\cdots s_{\b_{m-1}}s_{\a_{n-1}}
\end{array}
$$

\medskip
We just point out that
$$
T_x=
\begin{cases}
Z(G)
 & \textrm{for $x\in Z_m \cap wB$}\\
 (T^w)^\circ\times \<{\prod_{i=1}^{\ell}h_{\a_{2i-1}}(-1)}
 & \textrm{for $x\in Z_\ell\cap wB$, $\ell=1,\ldots m-1$ }
\end{cases}
$$
We get
\begin{center}
\vskip-20pt
$$
\begin{array}{|c||c|c|}
\hline
\O  & \l(\O)& \l(\hat\O) \\ 
\hline
\hline 
\begin{array}{c}
X_\ell\\
\ell=1,\ldots,m-1
\end{array}
& \displaystyle\sum_{i=1}^\ell n_{2i}\,\omega_{2i} &\\
\hline
\begin{array}{c}
X_m
\end{array}
& \displaystyle\sum_{i=1}^{m-1} n_{2i}\,\omega_{2i}+2n_n\omega_n& \displaystyle\sum_{i=1}^{m-1} n_{2i}\,\omega_{2i}+n_n\omega_n\\
\hline
\begin{array}{c}
X_m'
\end{array}
& \displaystyle\sum_{i=1}^{m-1} n_{2i}\,\omega_{2i}+2n_{n-1}\omega_{n-1}& \displaystyle\sum_{i=1}^{m-1} n_{2i}\,\omega_{2i}+n_{n-1}\omega_{n-1}\\
\hline
\begin{array}{c}
Z_\ell\\
\ell=1,\ldots,m-1
\end{array}
& \displaystyle\sum_{i=1}^{2\ell} n_i\omega_i,\ \sum_{i=1}^{\ell}n_{2i-1}\ \text{even}& \displaystyle\sum_{i=1}^{2\ell} n_i\omega_i \\
\hline
\begin{array}{c}
Z_m
\end{array}
&  \displaystyle\sum_{i=1}^n n_i\omega_i,\ \sum_{i=1}^{m}n_{2i-1}\ \text{even,\ $n_{n-1}+n_n$ even}& \displaystyle \sum_{i=1}^{n} n_i\omega_i\\
\hline
\end{array}
$$
\end{center}
\begin{center} Table \totable: $\l(\O)$, $\l(\hat\O)$ for unipotent classes in $D_n$, $n=2m$.
\end{center}
 In particular $\hat Z_m$ is a model homogeneus space, and in fact the principal one, by \cite{Luna}, 3.3 (4).

\medskip

\subsubsection{Semisimple classes in $D_n$, $n$ even $n=2m$}
Following the notation in \cite{CCC}, Tables 1, 5 we have
$$
\begin{array}{lllll}
D_\ell D_{n-\ell}&\corr &J_\ell,\quad \ell=1,\ldots,m&\corr &s_{\b_1}s_{\d_1}\cdots s_{\b_\ell}s_{\d_\ell}\\
T_1 A_{n-1}&\corr &K_m,\quad &\corr & s_{\b_1}\cdots s_{\b_m}\\
(T_1 A_{n-1})'&\corr & \{1,3,\ldots,n-3,n\}&\corr & s_{\b_1}\cdots s_{\b_{m-1}}s_{\a_{n-1}}
\end{array}
$$
There are two families of classes of semisimple elements with centralizer of type $T_1 A_{n-1}$: to distinguish them we wrote $T_1 A_{n-1}$ and $(T_1 A_{n-1})'$.
We get

\begin{center}
\vskip-20pt
$$
\begin{array}{|c|c||c|}
\hline
\O  & H & \l(\O) \\ 
\hline
\hline
\begin{array}{c}
\exp(\zeta\check\omega_1)\\
\zeta\in\C\setminus 2\pi i\Z
\end{array}
& T_1 D_{n-1}&2n_1\omega_1+n_2\omega_2\\
\hline
\begin{array}{c}
\exp(\pi i \check\omega_\ell)\\
\ell=2,\ldots,m-1
\end{array}
&D_\ell D_{n-\ell}&  \displaystyle\sum_{i=1}^{2\ell-1}2 n_i\omega_i+n_{2\ell}\omega_{2\ell}
\\
\hline
\begin{array}{c}
\exp(\pi i \check\omega_m)
\end{array}
& D_m D_{m}&  \displaystyle\sum_{i=1}^{n}2 n_i\omega_i\\
\hline
\begin{array}{c}
\exp(\zeta\check\omega_n)\\
\zeta\in\C\setminus 2\pi i\Z
\end{array}
& T_1 A_{n-1}&  \displaystyle \sum_{i=1}^{m-1} n_{2i}\,\omega_{2i}+2n_n\omega_n\\
\hline
\begin{array}{c}
\exp(\zeta\check\omega_{n-1})\\
\zeta\in\C\setminus 2\pi i\Z
\end{array}
&(T_1 A_{n-1})'&  \displaystyle\sum_{i=1}^{m-1} n_{2i}\,\omega_{2i}+2n_{n-1}\omega_{n-1}\\
\hline
\end{array}
$$
\end{center}
\begin{center} Table \totable: $\l(\O)$ for semisimple classes in $D_n$, $n=2m$.
\end{center}

\medskip

\subsubsection{Unipotent classes in $D_n$, $n$ odd, $n=2m+1$.}

\noindent
The center of $G$ is $\<{(\prod_{j=1}^{m}h_{\a_{2j-1}}(-1))h_{\a_{n-1}}(i)h_{\a_n}(-i)}$.
 From \cite{CCC} we have
$$
\begin{array}{llllll}
Z_\ell,&\ell=1,\ldots,m&\corr &J_\ell&\corr &s_{\b_1}s_{\d_1}\cdots s_{\b_\ell}s_{\d_\ell}\\
X_\ell,& \ell=1,\ldots,m &\corr &K_\ell&\corr & s_{\b_1}\cdots s_{\b_\ell}
\end{array}
$$

\medskip
We get
\begin{center}
\vskip-20pt
$$
\begin{array}{|c||c|c|}
\hline
\O  & \l(\O)& \l(\hat\O) \\ 
\hline
\hline 
\begin{array}{c}
X_\ell\\
\ell=1,\ldots,m-1
\end{array}
& \displaystyle\sum_{i=1}^\ell n_{2i}\,\omega_{2i} & \\
\hline
\begin{array}{c}
X_m
\end{array}
& \displaystyle\sum_{i=1}^{m-1} n_{2i}\,\omega_{2i}+n_{n-1}(\omega_{n-1}+\omega_n)& \\
\hline
\begin{array}{c}
Z_\ell\\
\ell=1,\ldots,m-1
\end{array}
& \displaystyle\sum_{i=1}^{2\ell} n_i\omega_i,\ \sum_{i=1}^{\ell}n_{2i-1}\ \text{even}& \displaystyle\sum_{i=1}^{2\ell} n_i\omega_i \\
\hline
\begin{array}{c}
Z_m
\end{array}
&  \displaystyle\sum_{i=1}^{n-2} n_i\omega_i+n_{n-1}(\omega_{n-1}+\omega_n),\ \sum_{i=1}^{m}n_{2i-1}\ \text{even}& \displaystyle \sum_{i=1}^{n-2} n_i\omega_i+n_{n-1}(\omega_{n-1}+\omega_n)
\\
\hline
\end{array}
$$
\end{center}
\begin{center} Table \totable:  $\l(\O)$, $\l(\hat\O)$ for unipotent classes in $D_n$, $n=2m+1$.
\end{center} 

\medskip

\subsubsection{Semisimple classes in $D_n$, $n$ odd, $n=2m+1$}
%Following the notation in \cite{CCC}, Tables 1, 5 we get
$$
\begin{array}{lllll}
D_\ell D_{n-\ell},\quad \ell=1,\ldots,m&\corr &J_\ell&\corr &s_{\b_1}s_{\d_1}\cdots s_{\b_\ell}s_{\d_\ell}\\
T_1 A_{n-1}&\corr &K_m&\corr & s_{\b_1}\cdots s_{\b_m}\\
\end{array}
$$
\medskip
We obtain

\begin{center}
\vskip-20pt
$$
\begin{array}{|c|c||c|}
\hline
\O  & H & \l(\O) \\ 
\hline
\hline
\begin{array}{c}
\exp(\zeta\check\omega_1)\\
\zeta\in\C\setminus 2\pi i\Z
\end{array}
& T_1 D_{n-1}&2n_1\omega_1+n_2\omega_2\\
\hline
\begin{array}{c}
\exp(\pi i \check\omega_\ell)\\
\ell=2,\ldots,m-1
\end{array}
&D_\ell D_{n-\ell}&  \displaystyle\sum_{i=1}^{2\ell-1}2 n_i\omega_i+n_{2\ell}\omega_{2\ell}
\\
\hline
\begin{array}{c}
\exp(\pi i \check\omega_m)
\end{array}
& D_m D_{m+1}&  \displaystyle \sum_{i=1}^{n-2}2 n_i\omega_i+n_{n-1}(\omega_{n-1}+\omega_n)\\
\hline
\begin{array}{c}
\exp(\zeta\check\omega_n)\\
\zeta\in\C\setminus 2\pi i\Z
\end{array}
&T_1 A_{n-1}&  \displaystyle\sum_{i=1}^{m-1} n_{2i}\,\omega_{2i}+n_{n-1}(\omega_{n-1}+\omega_n)\\
\hline
\end{array}
$$
\end{center}
\begin{center} Table \totable: $\l(\O)$ for semisimple classes in $D_n$, $n=2m+1$.
\end{center}

\subsection{Type $B_n$, $n\geq 2$.}

We put $m=[\frac{n}2]$. The center of $G$ is
 $
\<{h_{\a_{n}}(-1)}
 $.
We have $\omega_i=e_1+\cdots+e_i$ for $i=1,\ldots, n-1$,  $\omega_n=\frac12(e_1+\cdots+e_{n})$. We put $\g_\ell=e_{\ell}$,  $M_\ell=\{\ell+1,\ldots,n\}$ for $\ell=1,\ldots, n$ and  $\b_i=e_{2i-1}+e_{2i}$, $\d_i=e_{2i-1}-e_{2i}$, $J_\ell=\{2\ell+1,\ldots,n\}$,
$K_\ell=J_\ell\cup\{1,3,\ldots,2\ell-1\}$ for $i=1,\ldots, m$. 

\subsubsection{Unipotent classes in $B_n$, $n$ {even}, $n=2m$.} 

$$
\begin{array}{llllll}
Z_\ell,&\ell=1,\ldots,m&\corr &J_\ell&\corr &s_{\b_1}s_{\d_1}\cdots s_{\b_\ell}s_{\d_\ell}\\
X_\ell,&\ell=1,\ldots,m &\corr &K_\ell&\corr & s_{\b_1}\cdots s_{\b_\ell}
\end{array}
$$

We obtain
\medskip
\begin{center}
\vskip-20pt
$$
\begin{array}{|c||c|c|}
\hline
\O  & \l(\O)& \l(\hat\O) \\ 
\hline
\hline 
\begin{array}{c}
X_\ell\\
\ell=1,\ldots,m-1
\end{array}
& \displaystyle\sum_{i=1}^\ell n_{2i}\,\omega_{2i}&\\
\hline
\begin{array}{c}
X_m
\end{array}
& \displaystyle\sum_{i=1}^{m-1} n_{2i}\,\omega_{2i}+2n_n\omega_n& \displaystyle\sum_{i=1}^{m} n_{2i}\,\omega_{2i}\\
\hline
\begin{array}{c}
Z_\ell\\
\ell=1,\ldots,m-1
\end{array}
& \displaystyle\sum_{i=1}^{2\ell} n_i\omega_i,\ \sum_{i=1}^{\ell}n_{2i-1}\ \text{even}& \displaystyle\sum_{i=1}^{2\ell} n_i\omega_i \\
\hline
\begin{array}{c}
Z_m
\end{array}
&  \displaystyle\sum_{i=1}^n n_i\omega_i,\ \sum_{i=1}^{m}n_{2i-1}\ \text{even,\ $n_n$ even}& \displaystyle\sum_{i=1}^n n_i\omega_i,\ \text{ $n_n$ even}\\
\hline
\end{array}
$$
\end{center}
\begin{center} Table \totable: $\l(\O)$, $\l(\hat\O)$ for unipotent classes in $B_n$, $n=2m$.
\end{center}
 
\medskip

\subsubsection{Semisimple classes in $B_n$, $n$ even $n=2m$}
Following the notation in \cite{CCC}, Tables 1, 5 we have
$$
\begin{array}{llllll}
D_\ell B_{n-\ell},&\ell=1,\ldots,m&\corr &J_\ell&\corr &s_{\b_1}s_{\d_1}\cdots s_{\b_\ell}s_{\d_\ell}\\
D_\ell B_{n-\ell},& \ell=m+1,\ldots,n&\corr &M_{2(n-\ell)+1}&\corr &s_{\g_1}s_{\g_2}\cdots s_{\g_{2(n-\ell)+1}}\\
T_1 A_{n-1}&&\corr &\vuoto&\corr & w_0\end{array}
$$

We obtain
\medskip
\begin{center}
\vskip-20pt
$$
\begin{array}{|c|c||c|}
\hline
\O  & H & \l(\O) \\ 
\hline
\hline
\begin{array}{c}
\exp(\zeta\check\omega_1)\\
\zeta\in\C\setminus 2\pi i\Z, m\geq 2
\end{array}
& T_1 B_{n-1}&2n_1\omega_1+n_2\omega_2\\
\hline
\begin{array}{c}
\exp(\zeta\check\omega_1)\\
\zeta\in\C\setminus 2\pi i\Z, m=1
\end{array}
& T_1 B_{1}&2n_1\omega_1+2n_2\omega_2\\
\hline
\begin{array}{c}
\exp(\pi i \check\omega_\ell)\\
\ell=2,\ldots,m-1
\end{array}
&D_\ell B_{n-\ell}&  \displaystyle\sum_{i=1}^{2\ell-1}2 n_i\omega_i+n_{2\ell}\omega_{2\ell}
\\
\hline
\begin{array}{c}
\exp(\pi i \check\omega_m)\\
\end{array}
& D_m B_{m}&\displaystyle\sum_{i=1}^{n}2 n_i\omega_i\\
\hline
\begin{array}{c}
\exp(\pi i \check\omega_\ell)\\
\ell=m+1,\ldots,n
\end{array}
& D_\ell B_{n-\ell}&  \displaystyle\sum_{i=1}^{2(n-\ell)} 2n_i\omega_i+n_{2(n-\ell)+1}\omega_{2(n-\ell)+1}\\
\hline
\begin{array}{c}
\exp(\zeta\check\omega_n)\\
\zeta\in\C\setminus \pi i\Z
\end{array}
&T_1 A_{n-1}&  \displaystyle\sum_{i=1}^n n_i\omega_i,\  \text{\ $n_n$ even}\\
\hline
\end{array}
$$
\end{center}
\begin{center} Table \totable: $\l(\O)$ for semisimple classes in $B_n$, $n=2m$.
\end{center}

\medskip
\subsubsection{Mixed classes in $B_n$, $n$ even, $n=2m$}
From \cite{CCC}, Table 4, we get
$$
\begin{array}{llllll}
\s_n x_{\b_1}(1)\cdots x_{\b_m}(1)&& \corr & \vuoto& \corr & w_0\\
\s_n x_{\b_1}(1)\cdots x_{\b_\ell}(1),&\ell=1,\ldots,m-1& \corr & M_{2\ell+1}& \corr & s_{\g_1}\cdots s_{\g_{2\ell+1}}
\end{array}
$$
We obtain
\begin{center}
\vskip-20pt
$$
\begin{array}{|c||c|c|}
\hline
\O  & \l(\O)& \l(\hat\O) \\ 
\hline
\hline 
\begin{array}{c}
\s_n x_{\b_1}(1)\cdots x_{\b_\ell}(1)\\
\ell=1,\cdots,m-1
\end{array}
&\displaystyle\sum_{i=1}^{2\ell+1} n_i\omega_i
&\\
\hline
\begin{array}{c}
\s_n x_{\b_1}(1)\cdots x_{\b_m}(1)
\end{array}
&  \displaystyle \sum_{i=1}^n n_i\omega_i,\  \text{\ $n_n$ even}& \quad\displaystyle \sum_{i=1}^n n_i\omega_i\quad
\\
\hline
\end{array}
$$
\end{center}
\begin{center} Table \totable: $\l(\O)$, $\l(\hat\O)$ for mixed classes in $B_n$, $n=2m$.
\end{center}
In particular $\hat\O_{\s_n x_{\b_1}(1)\cdots x_{\b_m}(1)}$ is a model homogeneus space, and in fact the principal one, by \cite{Luna}, 3.3 (2).

\subsubsection{Unipotent classes in $B_n$, $n$ odd, $n=2m+1$.}

$$
\begin{array}{lllll}
Z_\ell&\corr &J_\ell,\quad \ell=1,\ldots,m&\corr &s_{\b_1}s_{\d_1}\cdots s_{\b_\ell}s_{\d_\ell}\\
Z_{m+1}&\corr &\vuoto &\corr&w_0=s_{\b_1}s_{\d_1}\cdots s_{\b_m}s_{\d_m}s_{\a_n}\\
X_\ell &\corr &K_\ell,\quad \ell=1,\ldots,m&\corr & s_{\b_1}\cdots s_{\b_\ell}
\end{array}
$$

\begin{lemma}\label{toni1bd} Let $w=s_{\b_1}\cdots s_{\b_\ell}
$
for $\ell=1,\ldots,m$. Then $T^w$ is connected. 
\end{lemma}
\pf For $\ell=1,\ldots,m$ we have $(1-w)P=\Z\<{\b_1,\ldots,\b_\ell}= \Z\<{\omega_{2i}\mid i=1,\ldots,\ell}$.\cvd

\medskip

\begin{lemma}\label{toni3bd}Let $w=s_{\b_1}\cdots s_{\b_\ell}s_{\d_1}\cdots s_{\d_\ell}$ for $\ell=1,\ldots,m
$. Then
 $$
 T^w=(T^w)^\circ\times \<{h_{\a_1}(-1)}\times\cdots\times  \<{h_{\a_{2\ell-1}}(-1)}
 $$ 
 \end{lemma}
 \pf For $\ell=1,\ldots,m$ we have 
$
(1-w)P=\Z\<{2\omega_1,\ldots,2\omega_{2\ell-1},\omega_{2\ell}}
$.\cvd

\medskip

For $\ell=1$ we get $T^w=\<{h_{\a_1}(-1)}\times (T^w)^\circ$.
In \cite{CCC} we exhibit the element  $x_{-\b_1}(1)
 x_{-\d_1}(1)\in \O\cap BwB\cap B^-$. We may therefore choose
$
x=n_{\b_1} n_{\d_1}h\,x_{\b_1}(2) x_{\d_1}(2)
$
for a certain $h\in T$. Then $h_{\a_1}(-1)\in C(x)$, so that $T_x=T^w$. 

Next we consider $Z_{m+1}$. We claim that $T_x=Z(G)$. Suppose for a contradiction that there is an involution $\sigma\in T_x\setminus Z(G)$. Then $x\in K=C(\sigma)$, and $K$ is the almost direct product $K_1 K_2$, of type $D_k B_{n-k}$, for some $k=1,\ldots,n$. We get an orthogonal decomposition $E=E_1\oplus E_2$ and a decomposition $x=x_1x_2\in K_1 K_2$. Then $-1=w_0=(w_1,w_2)$, where $w_i$ is the element of the Weyl group of $K_i$ corresponding to $x_i$ (the class of $x_i$ in $K_i$ is spherical). It follows that each $w_i=-1$, and $k$ is even. Then $x_1$ is in the class $Z_{k/2}$ of $K_1$ and $x_2$ in the class $Z_{m+1-k/2}$ of $K_2$. However, the product $x_1x_2$ is not in the class $Z_{m+1}$ of $G$ (since in $x_1x_2$ there are two rows with 3 boxes), a contradiction. Hence $T_x=Z(G)$. 

We now deal with $Z_\ell$, $\ell=2,\ldots,m$.
Here $\Psi_J$ has basis $\{\a_1,\ldots,\a_{2\ell-1},\g_{2\ell}\}$, and
 $C((T^w)^\circ)'$ is of type $B_{2\ell}$. From the construction in \cite{CCC}, proof of Theorem 2.11, we can find $x$ in the $D_{2\ell}$-subgroup $K$ of $C((T^w)^\circ)'$ generated by the long roots, that is the $D_{2\ell}$-subgroup with basis $\{\a_1,\ldots,\a_{2\ell-1},\b_\ell\}$.
We have 
$$
Z(K)=Z(G)\times 
\<\sigma\quad,\quad \sigma= \prod_{i=1}^{\ell}h_{\a_{2i-1}}(-1)
$$
By Lemma \ref{toni3bd},  $T_x=(T^w)^\circ\times (T_x\cap R)$,
where
$R= \<{h_{\a_1}(-1)}\times\cdots\times  \<{h_{\a_{2\ell-1}}(-1)}\leq K
$. Since $x$ lies in the maximal spherical unipotent class of $D_{2\ell}$, 
from the result obtained for this class, we have
$
T_x\cap R=R\cap Z(K)= \<\sigma
$,  hence  
$
 T_x=(T^w)^\circ\times \<\sigma
  $.
We have proved

\begin{proposition}\label{toni4bd} For $\ell=1,\ldots,m$ we have
{\rm
$$
\l(Z_\ell)=\left\{\sum_{i=1}^{2\ell} n_i\omega_i
\mid n_k\in \N,\ \sum_{i=1}^{\ell}n_{2i-1}\ \text{even}\right\}
$$}
Moreover
{\rm
$$
\l(Z_{m+1})=\left\{\sum_{i=1}^n n_i\omega_i
\mid n_k\in \N,\  \text{\ $n_n$ even}\right\}
$$}
\end{proposition}

\medskip
For the simply-connected cover we obtain
\begin{proposition}\label{toni4bduniversal} For $\ell=1,\ldots,m$ we have
{\rm
$$
\l(\hat Z_\ell)=\left\{\sum_{i=1}^{2\ell} n_i\omega_i
\mid n_k\in \N\right\}
$$}
Moreover
{\rm
$$
\l(\hat Z_{m+1})=\left\{\sum_{i=1}^n n_i\omega_i
\mid n_k\in \N\right\}
$$}
\end{proposition}
\pf 
Let $u\in Z_\ell$, with $\ell=1,\ldots,m+1$. If $C(u)^\circ =RC$ with $R=R_u(C(u))$,  $C$ connected reductive, then $C$ is of type 
$C_{\ell-1} D_{n-2\ell+1}$ (\cite{Carter2}, \S 13.1). In particular $C$ is semisimple since $n-2\ell+1$ is even. Hence $\l(\hat Z_\ell)$ is free by Lemma \ref{free}.

For $\ell=m+1$, we have $Z(G)\not\leq C(x)^\circ$. In fact, we can take $u=x_{\a_1}(1)x_{\a_3}(1)\cdots x_{\a_n}(1)$ in $Z_{m+1}$. Then $S=H_{\check\omega_2}H_{\check\omega_4}\cdots H_{\check\omega_{n-1}}$ is a maximal torus of $C(u)^\circ$, and since $Z(G)\cap S=\{1\}$, we get $C(u)=C(u)^\circ \times Z(G)$ by Lemma \ref{esclude}. We are left to deal with $\ell=1$. However for each $\ell$, the image $Q$ of $(u-1)^2$ in $V(\omega_1)$ (which is the natural module for $B_n$) has dimension 1, so $C(u)^\circ$ acts trivially on $Q$ by Lemma \ref{dim 1}, and  
$\omega_1\in \l(\hat Z_{\ell})$.\cvd

We summarize the results obtained in
\begin{center}
\vskip-20pt
$$
\begin{array}{|c||c|c|}
\hline
\O  & \l(\O)& \l(\hat\O) \\ 
\hline 
\begin{array}{c}
X_\ell\\
\ell=1,\ldots,m
\end{array}
& \displaystyle\sum_{i=1}^\ell n_{2i}\,\omega_{2i} &\\
\hline
\begin{array}{c}
Z_\ell\\
\ell=1,\ldots,m
\end{array}
& \displaystyle\sum_{i=1}^{2\ell} n_i\omega_i,\ \sum_{i=1}^{\ell}n_{2i-1}\ \text{even}&\quad \displaystyle\sum_{i=1}^{2\ell} n_i\omega_i \quad\\
\hline
\begin{array}{c}
Z_{m+1}
\end{array}
&  \sum_{i=1}^n n_i\omega_i,\  \text{\ $n_n$ even}& \displaystyle \sum_{i=1}^n n_i\omega_i\\
\hline
\end{array}
$$
\end{center}
\begin{center} Table \totable: $\l(\O)$, $\l(\hat\O)$ for unipotent classes in $B_n$, $n=2m+1$.
\end{center}
In particular $\hat Z_{m+1}$ is a model homogeneus space, and in fact the principal one, by \cite{Luna}, 3.3 (2).

\medskip
In section \ref{non-normal}, we shall determine the decomposition of the coordinate ring of the closure of $Z_{m+1}$. For this purpose we shall use the fact that if $x\in Z_{m+1}\cap w_0B$, then $\a_{n-1}$ occurs in $x$ (see the discussion before Proposition \ref{funzioni}). This can be checked by using the representative of $Z_{m+1}$ in $SO(2n+1)$ given 
in \cite{CCC}, proof of Theorem 12. 

\medskip

\subsubsection{Semisimple classes in $B_n$, $n$ odd $n=2m+1$}
Following the notation in \cite{CCC}, Tables 1, 5 we get
$$
\begin{array}{llllll}
D_\ell B_{n-\ell},& \ell=1,\ldots,m&\corr &J_\ell&\corr &s_{\b_1}s_{\d_1}\cdots s_{\b_\ell}s_{\d_\ell}\\
D_\ell B_{n-\ell},& \ell=m+1,\ldots,n&\corr &M_{2(n-\ell)+1}&\corr &s_{\g_1}s_{\g_2}\cdots s_{\g_{2(n-\ell)+1}}\\
T_1 A_{n-1}&&\corr &\vuoto&\corr & w_0
\end{array}
$$
and we obtain
\begin{center}
\vskip-20pt
$$
\begin{array}{|c|c||c|}
\hline
\O  & H & \l(\O) \\ 
\hline
\hline
\begin{array}{c}
\exp(\zeta\check\omega_1)\\
\zeta\in\C\setminus 2\pi i\Z
\end{array}
& T_1 B_{n-1}&2n_1\omega_1+n_2\omega_2\\
\hline
\begin{array}{c}
\exp(\pi i \check\omega_\ell)\\
\ell=2,\ldots,m
\end{array}
&D_\ell B_{n-\ell}&  \displaystyle\sum_{i=1}^{2\ell-1}2 n_i\omega_i+n_{2\ell}\omega_{2\ell}
\\
\hline
\begin{array}{c}
\exp(\pi i \check\omega_\ell)\\
\ell=m+2,\ldots,n
\end{array}
& D_\ell B_{n-\ell}&  \displaystyle\sum_{i=1}^{2(n-\ell)} 2n_i\omega_i+n_{2(n-\ell)+1}\omega_{2(n-\ell)+1}\\
\hline
\begin{array}{c}
\exp(\pi i \check\omega_{m+1})
\end{array}
& D_{m+1} B_{m}&  \displaystyle\sum_{i=1}^{n} 2n_i\omega_i\\
\hline
\begin{array}{c}
\exp(\zeta\check\omega_n)\\
\zeta\in\C\setminus \pi i\Z
\end{array}
&T_1 A_{n-1}&  \displaystyle\sum_{i=1}^n n_i\omega_i,\  \text{\ $n_n$ even}\\
\hline
\end{array}
$$
\end{center}
\begin{center} Table \totable: $\l(\O)$ for semisimple classes in $B_n$, $n=2m+1$.
\end{center}

\medskip
\subsubsection{Mixed classes in $B_n$, $n$ odd, $n=2m+1$}
From \cite{CCC}, Table 4, we get
$$
\begin{array}{cclll}
\s_n x_{\b_1}(1)\cdots x_{\b_\ell}(1),\ \ell=1,\ldots,m& \corr & M_{2\ell+1}& \corr & s_{\g_1}\cdots s_{\g_{2\ell+1}}
\end{array}
$$
and we obtain
\medskip
\begin{center}
\vskip-20pt
$$
\begin{array}{|c||c|}
\hline
\O  & \l(\O)= \l(\hat\O) \\ 
\hline
\hline 
\begin{array}{c}
\s_n x_{\b_1}(1)\cdots x_{\b_\ell}(1)\\
\ell=1,\cdots,m-1
\end{array}
&\quad \displaystyle\sum_{i=1}^{2\ell+1} n_i\omega_i\quad
 \\
\hline
\begin{array}{c}
\s_n x_{\b_1}(1)\cdots x_{\b_m}(1)
\end{array}
&\quad \displaystyle \sum_{i=1}^n n_i\omega_i,\  \text{\ $n_n$ even}\quad\\
\hline
\end{array}
$$
\end{center}
\begin{center} Table \totable: $\l(\O)$, $\l(\hat\O)$ for mixed classes in $B_n$, $n=2m+1$.
\end{center}

 \medskip
\subsection{Type $E_6$.}
We put
$$
\begin{array}{cclll}
\b_1 = (1,2,2,3,2,1),&
\b_2 =(1,0,1,1,1,1)\\
\b_3 =(0,0,1,1,1,0),&
\b_4 = (0,0,0,1,0,0)
\end{array}
$$

\subsubsection{Unipotent classes in $E_6$.}  

$$
\begin{array}{cclll}
A_1 & \corr & \{1,3,4,5,6\} & \corr & s_{\b_1}\\
2A_1& \corr & 
\{3,4,5\} & \corr & s_{\b_1} s_{\b_2}\\
3A_1& \corr & \vuoto\ & \corr & w_0=s_{\b_1}\cdots s_{\b_4}
\end{array}
$$

We obtain

\medskip
\begin{center}
\vskip-20pt
$$
\begin{array}{|c||c|}
\hline
\O  & \l(\O)= \l(\hat\O) \\ 
\hline
\hline 
\begin{array}{c}
A_1\end{array}
&\quad n_2\omega_2\quad
 \\
\hline
\begin{array}{c}
2A_1
\end{array}
&\quad n_1(\omega_1+\omega_6)+n_2\omega_2
\quad\\
\hline
\begin{array}{c}
3A_1
\end{array}
&\quad n_1(\omega_1+\omega_6)+n_3(\omega_3+\omega_5)+n_2\omega_2+n_4\omega_4
\quad\\
\hline
\end{array}
$$
\end{center}
\begin{center} Table \totable: $\l(\O)$, $\l(\hat\O)$ for unipotent classes in $E_6$.
\end{center}

\medskip
\subsubsection{Semisimple classes in $E_6$}

Following the notation in \cite{CCC}, Table 2, we have
$$
\begin{array}{lllllll}
A_1 A_5& \corr  & \vuoto\ & \corr & w_0\\
D_5\, T_1& \corr & 
\{3,4,5\} & \corr & s_{\b_1} s_{\b_2}\\
\end{array}
$$
We obtain
\medskip
\begin{center}
\vskip-20pt
$$
\begin{array}{|c|c||c|}
\hline
\O  & H & \l(\O) \\ 
\hline
\hline
\begin{array}{c}
\exp(\pi i \check\omega_2)
\end{array}
& A_1 A_5&n_1(\omega_1+\omega_6)+n_3(\omega_3+\omega_5)+2n_2\omega_2+2n_4\omega_4\\
\hline
\begin{array}{c}
\exp(\zeta\check\omega_1)\\
\zeta\in\C\setminus 2\pi i\Z
\end{array}
&D_5T_1&  n_1(\omega_1+\omega_6)+n_2\omega_2\\
\hline
\end{array}
$$
\end{center}
\begin{center} Table \totable: $\l(\O)$ for semisimple classes in $E_6$.
\end{center}
\medskip

\subsection{Type $E_7$.}

Here $Z(G)=\<{h_{\a_2}(-1)h_{\a_5}(-1)h_{\a_7}(-1)}$. We put
$$
\begin{array}{lll}
\b_1 = (2,2,3,4,3,2,1),\
\b_2 =(0,1,1,2,2,2,1),\
\b_3 =(0,1,1,2,1,0,0),\\
\b_4 = \a_7,\quad
\b_5=\a_5,\quad
\b_6=\a_3,\quad
\b_7=\a_2
\end{array}
$$

\subsubsection{Unipotent classes in $E_7$.}

$$
\begin{array}{cclll}
A_1 & \corr & \{2,3,4,5,6,7\} & \corr & s_{\b_1}\\
2A_1& \corr & \{2,3,4,5,7\} & \corr & s_{\b_1} s_{\b_2}\\
(3A_1)''& \corr & \{2,3,4,5\} & \corr & s_{\b_1} s_{\b_2}s_{\b_4}\\
(3A_1)'& \corr & \{2,5,7\} & \corr & s_{\b_1} s_{\b_2}s_{\b_3}s_{\b_6}\\
4A_1& \corr & \vuoto\ & \corr & w_0=s_{\b_1}\cdots s_{\b_7}
\end{array}
$$
We  obtain

\begin{center}
\vskip-20pt
$$
\begin{array}{|c||c|c|}
\hline
\O  & \l(\O)& \l(\hat\O) \\ 
\hline
\hline 
A_1& n_1\omega_1&  \\
\hline
\begin{array}{c}
\vspace{-0.4cm} 
2 A_1
\vspace{0.32cm} 
\end{array}
&  n_1\omega_1+n_6\omega_6&\\
\hline
\begin{array}{c}
\vspace{-0.4cm} 
(3A_1)''
\vspace{0.32cm} 
\end{array}
&  n_1\omega_1+n_6\omega_6+2n_7\omega_7
& n_1\omega_1+n_6\omega_6+n_7\omega_7\\
\hline
\begin{array}{c}
\vspace{-0.4cm} 
(3A_1)'
\vspace{0.32cm} 
\end{array}
&   n_1\omega_1+n_3\omega_3+n_4\omega_4+n_6\omega_6
&\\
\hline
\begin{array}{c}
\vspace{-0.4cm} 
4A_1
\vspace{0.32cm} 
\end{array}
&   \displaystyle\sum_{i=1}^7 n_i\omega_i,\  \text{$n_2+n_5+n_7$ even}&  \displaystyle\sum_{i=1}^7 n_i\omega_i\\
\hline
\end{array}
$$
\end{center}
\begin{center} Table \totable: $\l(\O)$, $\l(\hat\O)$ for unipotent classes in $E_7$.
\end{center}
In particular the simply-connected cover of ${4A}_1$ is a model homogeneus space, and in fact the principal one, by \cite{Luna}, 3.3 (8).

\begin{remark}\label{$E_7$ s.c.}{\rm From our description, it follows that $C(x)$ is connected for the classes $A_1$, $2A_1$ and $(3A_1)'$, while for $(3A_1)''$ and $4A_1$ we have $C(x)=C(x)^\circ\times Z(G)$.  This also follows from the tables in \cite{Ale}, where all unipotent classes are considered.}
\end{remark}

\subsubsection{Semisimple classes in $E_7$}

Following the notation in \cite{CCC}, Table 2, we have
$$
\begin{array}{cclll}
E_6T_1& \corr & \{2,3,4,5\} & \corr & s_{\b_1} s_{\b_2}s_{\b_4}\\
D_6A_1& \corr & \{2,5,7\} & \corr & s_{\b_1} s_{\b_2}s_{\b_3}s_{\b_6}\\
A_7& \corr & \vuoto\ & \corr & w_0
\end{array}
$$

We obtain
\medskip
\begin{center}
\vskip-20pt
$$
\begin{array}{|c|c||c|}
\hline
\O  & H & \l(\O) \\ 
\hline
\hline
\begin{array}{c}
\exp(\zeta\check\omega_7)\\
\zeta\in\C\setminus 2\pi i\Z
\end{array}
&E_6T_1&  n_1\omega_1+n_6\omega_6+2n_7\omega_7\\
\hline
\begin{array}{c}
\exp(\pi i \check\omega_1)
\end{array}
& D_6A_1&2n_1\omega_1+2n_3\omega_3+n_4\omega_4+n_6\omega_6\\
\hline
\begin{array}{c}
\exp(\pi i\check\omega_2)
\end{array}
& A_7&\displaystyle\sum_{i=1}^72n_i\omega_i\\
\hline
\end{array}
$$
\end{center}
\begin{center} Table \totable: $\l(\O)$ for semisimple classes in $E_7$.
\end{center}
\medskip

\subsection{Type $E_8$.}

We put
$$
\begin{array}{lll}
&\b_1 =   (2,3,4,6,5,4,3,2),\ 
\b_2=(2,2,3,4,3,2,1,0),\
\b_3=(0,1,1,2,2,2,1,0),\\
&\b_4=(0,1,1,2,1,0,0,0),\
\b_5=\a_7, \
\b_6=\a_5,\
\b_7=\a_3,\
\b_8=\a_2
\end{array}
$$

\subsubsection{Unipotent classes in $E_8$.}

$$
\begin{array}{cclll}
A_1 & \corr & \{1,2,3,4,5,6,7\} & \corr & s_{\b_1}\\
2A_1& \corr & \{2,3,4,5,6,7\} & \corr & s_{\b_1} s_{\b_2}\\
3A_1& \corr & \{2,3,4,5\} & \corr & s_{\b_1} s_{\b_2}s_{\b_3}s_{\b_5}\\
4 A_1& \corr & \vuoto\ & \corr & w_0=s_{\b_1}\cdots s_{\b_8}
\end{array}
$$
We have
$$
(1-w)P=
\begin{cases}
\Z\<{ \omega_8}
 & \textrm{for $w=s_{\b_1}$}\\
\Z\<{\omega_1,\omega_8}
 & \textrm{for $w=s_{\b_1}s_{\b_2}$}\\
\Z\<{\omega_1,\omega_6,2\omega_7,2\omega_8} & \textrm{for $w=s_{\b_1}s_{\b_2}s_{\b_3}s_{\b_5}$}
\end{cases}
$$
\medskip
\Class $3A_1$. Here $\Psi_J$ has basis $\{\a_7,\a_8,\b_2,\b_3\}$, 
 $K=C((T^w)^\circ)'$ is of type $D_{4}$ and has center
 $\<{h_{\a_3}(-1)h_{\a_5}(-1),h_{\a_2}(-1)h_{\a_3}(-1)}$ which is contained in $(T^w)^\circ$. Hence 
 $T_x=(T^w)^\circ$.

\medskip
\Class $4A_1$. 
We claim that $T_x=1$. Suppose for a contradiction there exists an involution $\sigma\in T_x$. Then $x\in K=C(\sigma)$. From the classification of involutions of $E_8$, it follows that $K$ is of type $D_8$ or $E_7 A_1$. The class of $x$ in $K$ is spherical, and by the uniqueness of Bruhat decomposition, $x$ lies over the longest element of the Weyl group of $K$, which is $w_0$. By comparison of weighted Dynkin diagrams, the spherical unipotent class of $K$ over $w_0$ does not correspond to the class $4A_1$ of $E_8$, a contradiction.

\medskip
We have shown that in all cases $T_x=(T^w)^\circ$, so that $C(x)$ is connected, as also follows from 
\cite{Carter2}, p. 405. We have

\begin{center}
\vskip-20pt
$$
\begin{array}{|c||c|}
\hline
\O  & \l(\O)= \l(\hat\O) \\ 
\hline
\hline 
A_1& n_8\omega_8 \\
\hline
\begin{array}{c}
\vspace{-0.4cm} 
2A_1
\vspace{0.32cm} 
\end{array}
&  n_1\omega_1+n_8\omega_8 \\
\hline
\begin{array}{c}
\vspace{-0.4cm} 
3A_1
\vspace{0.32cm} 
\end{array}
&  n_1\omega_1+n_6\omega_6+n_7\omega_7+n_8\omega_8\\
\hline
\begin{array}{c}
\vspace{-0.4cm} 
4A_1
\vspace{0.32cm} 
\end{array}
&   \displaystyle\sum_{i=1}^8 n_i\omega_i \\
\hline
\end{array}
$$
\end{center}
\begin{center} Table \totable: $\l(\O)$, $\l(\hat\O)$ for unipotent classes in $E_8$.
\end{center}
In particular $4A_1$ is a model homogeneus space (see \cite{vogan}, Theorem 1.1), and in fact the principal one, by \cite{Luna}, 3.3 (9).
\medskip

\medskip

\subsubsection{Semisimple classes in $E_8$.}

Following the notation in \cite{CCC}, Table 2, we have
$$
\begin{array}{cclll}
A_1E_7& \corr & \{2,3,4,5\} & \corr &s_{\b_1}s_{\b_2}s_{\b_3}s_{\b_5}\\
D_8& \corr & \vuoto\ & \corr & w_0
\end{array}
$$
We obtain
\medskip
\begin{center}
\vskip-20pt
$$
\begin{array}{|c|c||c|}
\hline
\O  & H & \l(\O) \\ 
\hline
\hline
\begin{array}{c} 
\exp(\pi i \check\omega_8)
\end{array}
&A_1E_7&n_1\omega_1+n_6\omega_6+2n_7\omega_7+2n_8\omega_8
\\
\hline
\begin{array}{c}
\exp(\pi i\check\omega_1)
\end{array}
& D_8&\displaystyle\sum_{i=1}^82n_i\omega_i\\
\hline
\end{array}
$$
\end{center}
\begin{center} Table \totable: $\l(\O)$ for semisimple classes in $E_8$.
\end{center}
\medskip

\subsection{Type $F_4$.}

We put
$$
\begin{array}{ll}
\b_1 =(2,3,4,2),&\b_2=(0,1,2,2),\\
\b_3=(0,1,2,0),&\b_4=(0,1,0,0)
\end{array}
$$

\subsubsection{Unipotent classes in $F_4$.}

$$
\begin{array}{cclll}
A_1 & \corr & \{2,3,4\} & \corr & s_{\b_1}\\
\tilde A_1& \corr & \{2,3\} & \corr & s_{\b_1} s_{\b_2}\\
A_1+\tilde A_1& \corr & \vuoto\ & \corr & w_0=s_{\b_1}\cdots s_{\b_4}
\end{array}
$$
We obtain
\begin{center}
\vskip-20pt
$$
\begin{array}{|c||c|c|}
\hline
\O  & \l(\O)& \l(\hat\O) \\ 
\hline
\hline 
A_1& n_1\omega_1&  \\
\hline
\begin{array}{c}
\vspace{-0.4cm} 
\tilde A_1
\vspace{0.32cm} 
\end{array}
&  n_1\omega_1+2n_4\omega_4&\ n_1\omega_1+n_4\omega_4\ \\
\hline
\begin{array}{c}
\vspace{-0.4cm} 
A_1+\tilde A_1
\vspace{0.32cm} 
\end{array}
&   n_1\omega_1+n_2\omega_2+2n_3\omega_3+2n_4\omega_4
& \\
\hline
\end{array}
$$
\end{center}
\begin{center} Table \totable: $\l(\O)$, $\l(\hat\O)$ for unipotent classes in $F_4$.
\end{center}

\medskip

\subsubsection{Semisimple classes in $F_4$.}

Following the notation in \cite{CCC}, Table 2, we have
$$
\begin{array}{cclll}
A_1 C_3& \corr & \vuoto& \corr & w_0\\
B_4& \corr & \{1,2,3\} & \corr & s_{\gamma_1}
\end{array}
$$
where $\gamma_1$ is the highest short root $(1,2,3,2)$.

We obtain
\medskip
\begin{center}
\vskip-20pt
$$
\begin{array}{|c|c||c|}
\hline
\O  & H & \l(\O) \\ 
\hline
\hline
\begin{array}{c}
\exp(\pi i \check\omega_1)
\end{array}
& A_1C_3&\quad\displaystyle\sum_{i=1}^4 2n_i\omega_i\quad\\
\hline
\begin{array}{c}
\exp(\pi i \check\omega_4)
\end{array}
& B_4&n_4\omega_4\\
\hline
\end{array}
$$
\end{center}
\begin{center} Table \totable: $\l(\O)$ for semisimple classes in $F_4$.
\end{center}
\medskip

\subsubsection{Mixed class in $F_4$.}

We put $f_2=\exp(\pi i \check\omega_4)$. Then following \cite{CCC}, Table 4
$$
\begin{array}{cclll}
\O_{f_2x_{\b_1}(1)}& \corr & \vuoto& \corr & w_0
\end{array}
$$
Assuming the existence of an involution in $T_x$ we get a contradiction, proving therefore that $T_x=1$. Hence
\begin{center}
\vskip-20pt
$$
\begin{array}{|c|c|}
\hline
\O  & \l(\O)=\l(\hat\O) \\ 
\hline
\hline
\begin{array}{c}
f_2x_{\b_1}(1)
\end{array}
&\displaystyle\sum_{i=1}^4 n_i\omega_i\\
\hline
\end{array}
$$
\end{center}
\begin{center} Table \totable: $\l(\O)$ for the mixed class in $F_4$.
\end{center}
\medskip
In particular $\O_{f_2x_{\b_1}(1)}$ is a model homogeneus space, and in fact the principal one, by \cite{Luna}, 3.3 (6), see also \cite{Luna} p. 300.

\medskip

\subsection{Type $G_2$.}
We put 
$
\b_1 =(3,2),\
\b_2=\a_1 
$.

\subsubsection{Unipotent classes in $G_2$.}

$$
\begin{array}{cclll}
A_1 & \corr & \{1\} & \corr &s_{\b_1}\\
\tilde A_1& \corr & \vuoto\ & \corr & w_0=s_{\b_1} s_{\b_2}
\end{array}
$$
We get
\begin{center}
\vskip-20pt
$$
\begin{array}{|c||c|}
\hline
\O  & \l(\O)= \l(\hat\O) \\ 
\hline
\hline 
A_1& n_2\omega_2 \\
\hline
\begin{array}{c}
\vspace{-0.4cm} 
\tilde A_1
\vspace{0.32cm} 
\end{array}
&\quad   n_1\omega_1+n_2\omega_2\quad \\
\hline
\end{array}
$$
\end{center}
\begin{center} Table \totable: $\l(\O)$, $\l(\hat\O)$ for unipotent classes in $G_2$.
\end{center}
In particular $\tilde A_1$ is a model homogeneus space, and in fact the principal one, by \cite{Luna}, 3.3 (5).

\medskip
Using the embedding of $G$ into $SO(7)$, one can determine explicitly an $x\in \O\cap w_0B$, where $\O= \tilde A_1$. Then one can check that both $\a_1$ and $\a_2$ occur in $x$ (see the discussion before Proposition \ref{funzioni}). This fact will be used in section \ref{non-normal} to determine $\C[\ov \O]$.

\subsubsection{Semisimple classes in $G_2$.}

Following the notation in \cite{CCC}, Table 2, we have
$$
\begin{array}{cclll}
A_1 \tilde A_1& \corr & \vuoto& \corr & w_0\\
A_2& \corr & \{2\} & \corr & s_{\gamma_1}
\end{array}
$$
where $\gamma_1$ is the highest short root $(2,1)$.

The group $G_2$ has 1 class of involutions. However there is also a class of elements of order 3 which is spherical.
We obtain
\medskip
\begin{center}
\vskip-20pt
$$
\begin{array}{|c|c||c|}
\hline
\O  & H & \l(\O) \\ 
\hline
\hline
\begin{array}{c}
\exp(\pi i \check\omega_2)
\end{array}
&A_1 \tilde A_1&\displaystyle\sum_{i=1}^2 2n_i\omega_i\\
\hline
\begin{array}{c}
\exp(\frac{2\pi i}3 \check\omega_1)
\end{array}
& A_2&n_1\omega_1\\
\hline
\end{array}
$$
\end{center}
\begin{center} Table \totable: $\l(\O)$ for semisimple classes in $G_2$.
\end{center}
\medskip

\section{The coordinate ring of $\ov\O$}\label{non-normal}

In this section we determine the decomposition of $\C[\ov\O]$ into simple $G$-modules, where $\ov\O$ is the closure of a spherical conjugacy class. Normality of conjugacy classes' closures has been deeply investigated.
For a survey on this topic, see \cite{Kraft}, \S 8, \cite{Broer2}, 7.9, Remark (iii).
The first observation is that the problem is reduced to unipotent conjugacy classes in $G$ (\cite{Kraft}, 8.1). In the following we are interested only in spherical conjugacy classes, and I  recall the facts in this context. It is known that the closure of the minimal nilpotent orbit is always normal (\cite{VP}, Theorem 2). Hesselink (\cite{Hes}) proved normality for several small orbits in the classical cases and certain orbits for the exceptional cases: namely, following the notation in \cite{Carter2},
$A_1$ and $2A_1$ in $E_6$, $A_1$, $2A_1$ and $(3A_1)''$ in $E_7$, $A_1$ and $2A_1$ in $E_8$, $A_1$ and $\tilde A_1$ in $F_4$, $A_1$ in $G_2$.

The classical groups have been considered in \cite{KP1}, \cite {KP2}: for the special linear groups the closure of every conjugacy class is normal. For the symplectic and orthogonal groups there exist conjugacy classes with non-normal closure. However every spherical conjugacy class in the symplectic group has normal closure, since from the classification we know that the unipotent spherical conjugacy classes have only 2 columns (see also  \cite{Hes}, \S 5, Criterion 2). For special orhogonal groups the results in \cite {KP2} left open the cases of the very even unipotent classes. E. Sommers proved that these have normal closure in \cite{Som2}. Taking into account the results in  \cite {KP2} and \cite{Som2} it follows that every unipotent spherical conjugacy class in type $D_n$ and $B_n$ has normal closure except for the maximal class $Z_{m+1}$ in $B_n$, when $n=2m+1$, $m\geq 1$. From this and the classification of spherical conjugacy classes, it follows that every spherical conjugacy class has normal closure, except for the above mentioned class in $B_{2m+1}$.

For the exceptional groups, besides the results on the minimal orbit and
Hesselink's results,  in \cite{LS} it is shown that the orbit $\tilde A_1$ in $G_2$ has a non-normal closure (see also \cite{Kraft}): here there is bijective normalization, contrary to the case of $Z_{m+1}$ in $B_{2m+1}$ where the closure is branched in codimension 2. In \cite{Broer1} the case of type $F_4$ is completely handled, and it follows that  every spherical conjugacy class has normal closure. The same holds for $E_6$, as follows from \cite{Som1} where every nilpotent orbit is considered. For the remaining nilpotent orbits in $E_7$ and $E_8$, in \cite{Broer2}, 7.9, Remark (iii), A. Broer gives a list of orbits with normal closure. Among these there are all spherical nilpotent orbits in $E_7$ and $E_8$. We may therefore state

\begin{theorem}\label{riassunto-normali} Let $\O$ be a spherical conjugacy class. Then $\ov\O$ is normal except for
the class $Z_{m+1}$ in $B_{2m+1}$ ($m\geq 1$) and the class $\tilde A_1$ in $G_2$.\cvd
\end{theorem}

\begin{remark}{\rm 
In \cite{CDCM}, Example 4.4, Proposition 4.5, the authors prove normal closure for nilpotent orbits of height 2.
}
\end{remark}

\begin{remark}{\rm 
In \cite{pany2}, 6.1, normality of $\cal N^{\rm sph}$ (the union of all spherical nilpotent orbits, which is in fact the closure of the unique maximal spherical nilpotent orbit) is discussed.
}
\end{remark}

\begin{remark}{\rm 
From (\ref{catena}) and Corollary \ref{funzioni18} it is possible to prove normality of $\ov\O$ in certain cases. For instance in type $C_n$ from Table 3 we get $\l(X_\ell)=2P^+_w$ for every unipotent class $X_\ell$. From (\ref{catena})  it follows that $\l(\ov\O)=\l(\O)$, so that $\ov\O$ is normal.
}
\end{remark}

We recall that in general $\C[\O]$ is the integral closure of $\C[\ov\O]$ in its field of fractions and that $\C[\ov\O]=\C[\O]$ if and only if $\ov\O$ is normal (\cite{Jan}, Proposition and Corollary in 8.3). By Theorem \ref{riassunto-normali}, to describe the decomposition of $\C[\ov\O]$ we are left to deal with $Z_{m+1}$ in $B_{2m+1}$ and with $\tilde A_1$ in $G_2$.
We use the notation and the tables from section 4 for the cases $B_{2m+1}$ and $G_2$.

\begin{theorem}\label{Non-normal} Let $\O=Z_{m+1}$ in $B_{n}$,  $n=2m+1$, $m\geq 1$. Then
{\rm
$$
\l(\ov\O)=\left\{\sum_{i=1}^{2m} n_i\omega_i\mid  \sum_{i=1}^{m}n_{2i-1}\ \text{even}\right\}\cup\left\{\sum_{i=1}^n n_i\omega_i\mid   \text{\ $n_n$ even, $n_n\geq 2$}\right\}
$$}
\end{theorem}
\pf Considering the ($G$-equivariant) restriction $r:\C[\ov\O]\to\C[\ov{Z_m}]=\C[Z_m]$, we get $\left\{\sum_{i=1}^{2m} n_i\omega_i
\mid  \sum_{i=1}^{m}n_{2i-1}\ \text{even}\right\}\leq \l(\ov\O)$. In particular for every even $j$, $\omega_j\in \l(\ov\O)$, and for every pair of odd $j$, $k$, with $1\leq j\leq k<n$, $\omega_j+\omega_k\in \l(\ov\O)$. By Corollary \ref{funzioni3}, we have $2\omega_n\in \l(\ov\O)$. We show that $\omega_j+2\omega_n\in \l(\ov\O)$ for every odd $j$, $j<n$.
We have
$
2\omega_{n-1}-\a_{n-1}=\omega_{n-2}+2\omega_n
$
and since $\a_{n-1}$ occurs in $x\in w_0B\cap\O$, by Corollary \ref{funzioni18}, we get $\omega_{n-2}+2\omega_n\in \l(\ov\O)$. Let $j$ be odd, $j<n-2$. Then 
$\omega_j+2\omega_n+2\omega_{n-2}\in\l(\ov\O)$ since
$\omega_{n-2}+2\omega_n$ and $\omega_j+\omega_{n-2}$ are in $\l(\ov\O)$.

There exists $B$-eigenvectors $F$, $H$ in $\C[\ov\O]$ of weights $\omega_j+2\omega_n+2\omega_{n-2}$, $2\omega_{n-2}$ respectively. Then $F/H$ is a rational function on $\ov\O$ of weight $\omega_j+2\omega_n$ defined at least on $\O$. However $2\omega_{n-2}$ is also a weight in $\l(Z_m)$, so that $H$ is non-zero on the dense $B$-orbit $\v$  in $Z_m$. Hence $F/H$ is defined on $\v$, and it is zero on $\v$, since $F$ is zero on $Z_m$, $\omega_j+2\omega_n+2\omega_{n-2}$ not being in $\l(Z_m)$. It follows that $F/H$ is defined on $Z_m$, so that it is a regular function on $\O\cup Z_m$. By \cite{KP2}, Theorem 16.2, (iii), $F/H$ extends to $\ov\O$, and  $\omega_j+2\omega_{n}$ lies in $ \l(\ov\O)$. We have shown that
$$
\l(\ov\O)\geq\left\{\sum_{i=1}^{2m} n_i\omega_i
\mid \sum_{i=1}^{m}n_{2i-1}\ \text{even}\right\}\cup\left\{\sum_{i=1}^n n_i\omega_i
\mid  \text{\ $n_n$ even, $n_n\geq 2$}\right\}
$$ 
We prove that also the opposite inclusion holds. 
Assume $\l=\sum_{i=1}^n n_i\omega_i\in \l(\ov\O)$. Since $\l(\ov\O)\leq \l(\O)$, we have $n_n$ even. If $n_n\not=0$ we are done. So assume $n_n=0$. Let $y\in Z_{m+1}\cap U^-\cap Bw_0B$.
We observe that $y_1:=\lim_{z\to 0}h_{\a_n}(z)^{-1}yh_{\a_n}(z)$ exists, and lies in $Z_m\cap U^-\cap BwB$, where $w=w(Z_m)$ (in \cite{CCC} we give representatives for both classes in $SO(2n+1)$, so that this may be checked directly). Now let $F:\ov\O\to \C$ be a highest weight vector of weight $\l$, with $F(y)=1$. Then $F(y_1)=1$, since $\l(h_{\a_n}(z))=1$ for every $z\in \C^\ast$. Since $x_1\in Z_m\cap wB$ lies in the $B$-orbit of $y_1$, we have $F(x_1)\not=0$. But $\s=\prod_{i=1}^{m}h_{\a_{2i-1}}(-1)\in C(x_1)$, so that
$F(x_1)=F(\s x_1\s)=\l(\s)F(x_1)$ implies $\l(\s)=1$, and we are done.
\cvd

\begin{theorem}\label{Non-normal2} Let $\O=\tilde A_1$ in $G_2$. Then
$
\l(\ov\O)
$
is the submonoid of $\l(\O)$ generated by $2\omega_1, 3\omega_1, \omega_2$.
\end{theorem}
\pf We know that $\omega_1\in \l(\O)$ and it follows from the proof of  \cite{LS}, Theorem 3.13, that $\omega_1\not\in\l(\ov\O)$.
We have
$$
2\omega_1-\a_1=\omega_2\quad,\quad
2\omega_2-\a_2=3\omega_1
$$
hence, by Corollary \ref{funzioni3} and \ref{funzioni18}, we get $2\omega_1$, $3\omega_1$, $\omega_2\in \l(\ov\O)$, since both $\a_1$, $\a_2$ occur in $x\in w_0B\cap\O$. Suppose for a contradiction that $\omega_1+n\omega_2\in \l(\ov\O)$ for a certain $n\in \N$. There exists $B$-eigenvectors $F$, $H$ in $\C[\ov\O]$ of weights $\omega_1+n\omega_2$, $n\omega_2$ respectively. Then $F/H$ is a rational function on $\ov\O$ of weight $\omega_1$ defined at least on $\O$. However $n\omega_2$ is also a weight in $\l(A_1)$, so that $H$ is non-zero on the dense $B$-orbit $\v$  in $A_1$. Hence $F/H$ is defined on $\v$, and it is zero on $\v$, since $F$ is zero on $A_1$, because $\omega_1+n\omega_2$ is not in $\l(A_1)$. It follows that $F/H$ is defined on $A_1$. But $A_1$ has normal closure, so that $F/H$ is defined on the closure of $A_1$, and then on $\ov\O$, so that there is in $\C[\ov\O]$ a $B$-eigenvector of weight $\omega_1$, a contradiction.  \cvd

\section{The general case}\label{generale}

Let $G$ be as usual simply-connected, $D\leq Z(G)$, $\ov G=G/D$, $\pi:G\to \ov G$ the canonical projection. For $g\in G$ we put $\ov g=\pi(g)$. We give a procedure to describe the coordinate ring of $\O_{\ov p}$, where $\O_{\ov p}$ is a spherical conjugacy class of $\bG$. Passing to $G$, we have to consider the quotient $G/\pi^{-1}(C_{\bG}(\ov p))$. Let $p=sv$ be the Jordan-Chevalley decomposition of $p$, $w=w(\O_p)$. We may assume $s\in T$. Let  $W_{s,D}=\{w\in W\mid wsw^{-1}=zs, z\in D\}$, and $N_{s,D}\leq N$ such that $N_{s,D}/T=W_{s,D}$. Then  $\pi^{-1}(C_{\bG}(\ov p))=C(v)\cap N_{s,D}C(s)$.
Reasoning as in \cite{131}, Corollary II, 4.4, we have a homomorphism $\pi^{-1}(C_{\bG}(\ov p))\to D$, $g\mapsto [g,p]$ with kernel $C(p)$. 

Let $y\in \O_p\cap BwB$ be such that $L=L_J$ is adapted to $C(y)$. If $H=\pi^{-1}(C_{\bG}(\ov y))$, then 
$\l(\O_{\ov p})=\l(G/H)=\{\l\in P^+_w\mid \l(T\cap H)=1\}$ by Corollary \ref{Gen. saturo}. Let $x\in \O_p\cap wB$, 
$x=\dot w u$, with $u\in U$ and let $T_{x,D}=T\cap \pi^{-1}(C_{\bG}(\ov x))$. 
By Proposition \ref{nane5}, we get $T\cap H=T_{x,D}$, hence
\begin{equation}
\l(\O_{\ov x})=\{\l\in P^+_w\mid \l(T_{x,D})=1\}
\label{non s.c.}
\end{equation}
Let $T_D^w=\{t\in T\mid wtw^{-1}=zt, z\in D\}$. From the Bruhat decomposition, we get $T_{x,D}\leq T^w_D$.
Moreover
since $w$ is an involution, for $t\in T^w_D$  we have $t=w^2tw^{-2}=z^2t$, so that $z^2=1$. In particular $\pi^{-1}(C_{\ov G}(\ov s))=N_{s,D_2}C(s)$,  $T_{D}^w= T_{D_2}^w$,
 where $D_2=D\cap T_2$.

Let $t\in T$ and write $t=ab$, with $a\in (T^w)^\circ$, $b\in (S^w)^\circ$. Then $wtw^{-1}=tz$ with $z\in D_2$ if and only if  $z=b^2$. Since $(S^w)^\circ$ is connected, we get $T_{D}^w=T_{D_2\cap (S^w)^\circ}^w$ and
$$
\frac{\pi^{-1}(C_{\bG}(\ov x))}{C(x)}\cong\frac{T_{x,D}}{T_x}\hookrightarrow \frac{T^w_{D}}{T^w}\cong D_2\cap (S^w)^\circ
$$
with $T_x=T^w\cap C(u)$, $T_{x,D}=T^w_D\cap C(u)$.
In particular, if $D_2\cap (S^w)^\circ=1$, then $\l(\O_{\ov x})=\l(\O_x)$.
This equality means that $x$ is not conjugate to $zx$ for any $z\in D_2$, $z\not=1$, and this may be directly checked in many cases, for instance in type $A_n$ or $C_n$ (and of course always holds for $x$ unipotent). However, to deal with orthogonal groups and $E_7$, we determined explicitly the cases when $D_2\cap (S^w)^\circ$ is non-trivial, and in each case we determined $T_{x,D}$ and therefore $\l(\O_{\ov x})$. 

Here we just observe that
if $D_2\cap (S^w)^\circ\not=1$, then $D_2\cap (S^w)^\circ\cong \Z/2\Z$, except possibly for $D=Z(G)$ in type $D_{n}$, $n=2m$. It turns out that in this case for $\exp(\pi i\check\omega_m)$, we have $T_x=T_2$ and  $T_{x,Z(G)}/T_x\cong \Z/2\Z\times \Z/2\Z$. More precisely 
$$
T_{x,Z(G)}=T_{Z(G)}^{w_0}=T_2\,\<{h_{\a_{n-1}}(i)h_{\a_{n}}(i),\prod_{i=1}^mh_{\a_{2i-1}}(i)}
$$
so that in $G/Z(G)=PSO(2n)$, $n=2m$, 
$$
\l(\O_{\ov{\exp(\pi i\check\omega_m)}})=\left\{\sum_{k=1}^{n}2m_k\omega_k\mid m_k\in \N,\  \text{\ $m_{n-1}+m_n$ and\, $\sum_{i=1}^{m} m_{2i-1}$ even}\right\}
$$
We add that for $SO(2n+1)$, $n\geq 1$ and $b_\l=\mbox{diag}(1,\l I_{n},\l^{-1}I_{n})$, $\l\not=\pm 1$, $\O_{b_\l}$ is a model orbit, and in fact the principal one by \cite{Luna}, 3.3 $(2')$.

We conclude by presenting the results for $E_7$. 

\subsection{Type $E_7$, $D=Z(G)$}

In this case $Z(G)=\<z$, where $z=h_{\a_2}(-1)h_{\a_5}(-1)h_{\a_7}(-1)=\exp(2\pi i\check\omega_2)=\exp(2\pi i\check\omega_7)$. 

There are 3 elements of the Weyl group to be considered and
only for $w=s_{\b_1} s_{\b_2}s_{\b_4}$ and $w=w_0$ we have $z\in (S^w)^\circ$. 

\medskip
\Class of type $A_7$, $w=w_0$. Here 
$
x=n_{\b_1}\cdots n_{\b_7}
$,
$$
T_{Z(G)}^{w_0}=T_2\, \<{\exp(\pi i \check\omega_2)}=
T_2\, \<{h_{\a_2}(i)h_{\a_5}(i)h_{\a_7}(i)}
$$
since $\exp(\pi i \check\omega_2)\in (S^{w_0})^\circ=T$ and $\exp(\pi i \check\omega_2)^2=z$.
\begin{proposition}\label{E7ad1} Let  $G$ be of type $E_7$, $D=Z(G)$, then
$$
\l(\O_{\ov{\exp(\pi i \check\omega_2)}}) =\left\{\sum_{i=1}^72n_i\omega_i\mid  \text{$n_2+n_5+n_7$ even}\right\}
$$
\end{proposition}
\pf This follows from the fact that $T_{x,Z(G)}= T^{w_0}_{Z(G)}$.
\cvd

\medskip
\Classes of type $E_6T_1$,  $w=s_{\b_1} s_{\b_2}s_{\b_4}$, $T^w=(T^w)^\circ\times \<{h_{\a_7}(-1)}=(T^w)^\circ\times Z(G)$.

We have 
$
T_{Z(G)}^w=T^w\<{\exp(\pi i\check\omega_7)}=T^w\<{h_{\a_1}(-1)h_{\a_7}(i)}
$.
If $\zeta\in\C\setminus 2\pi i\Z$, then
$$
x_\zeta=n_{\b_1}n_{\b_2}n_{\a_7}hx_{\b_1}(\xi)x_{\b_2}(\xi)x_{\a_7}(\xi)\in\O_{\exp(\zeta\check\omega_7)}\cap n_{\b_1}n_{\b_2}n_{\a_7}B
$$
for a certain $h\in T$, with $\xi=\frac{1+e^\zeta}{1-e^\zeta}$, so that
$$
T_{x_\zeta,Z(G)}=
\begin{cases}
T^w_{Z(G)}
 & \textrm{if $\zeta\in \pi i\Z\setminus 2\pi i\Z$}\\
T^w
 & \textrm{if $\zeta\in \C\setminus \pi i\Z$}
\end{cases}
$$
since $\a_7(\exp(\pi i\check\omega_7))=-1$.
\begin{proposition}\label{E7ad2} Let  $G$ be of type $E_7$, $D=Z(G)$, then
$$
\l(\O_{\ov{\exp(\zeta\check\omega_7)}})=
\begin{cases}
\left\{n_1\omega_1+n_6\omega_6+2n_7\omega_7\mid \text{$n_1+n_7$ even}\right\}
& \textrm{if $\zeta\in \pi i\Z\setminus 2\pi i\Z$}\\
\left\{n_1\omega_1+n_6\omega_6+2n_7\omega_7\right\}
 & \textrm{if $\zeta\in \C\setminus \pi i\Z$}
\end{cases}
$$
\end{proposition}
\cvd

\medskip
%\begin{addendum}\label{correzione}
\noindent
{\bf Addendum}
{\rm 
\ In \cite{CCC}, Remark 5, we stated that if $\pi_1:G\to G/U$ is the canonical projection, and $\O$ is a spherical conjugacy class, then ${\pi_1}_{|\O}:\O\to G/U$ has finite fibers. This is not correct, and one can only say that ${\pi_1}_{|\O}$ has generically finite fibers (if $w=w(\O)$, and $g\in \O\cap BwB$, then $\pi_1^{-1}(gU)$ has $\o{T^w/T_x}$ elements, where $x\in \O\cap wB$).
}
%\end{addendum}

\end{document}